\documentclass[12pt,twoside,leqno]{amsart}
\usepackage{amsmath,amscd,amssymb,amsfonts}

\usepackage[all]{xy}
\usepackage{graphicx}

\usepackage[latin1]{inputenc}
\usepackage[T1]{fontenc}

\reversemarginpar

\usepackage{amsmath}
\usepackage{amssymb}
\usepackage{latexsym}
\usepackage{amsfonts}
\usepackage{epsfig}
\usepackage{picinpar}
\usepackage[osf,sc]{mathpazo}
\usepackage{microtype}
\usepackage{textcomp}
\newcommand{\N}{\mathbb{N}}
\newcommand{\Z}{\mathbb{Z}}
\newcommand{\Q}{\mathbb{Q}}
\newcommand{\R}{\mathbb{R}}
\newcommand{\C}{\mathbb{C}}
\renewcommand{\H}{\mathbb{H}}
\renewcommand{\P}{\mathbb{P}}

\usepackage{graphicx}

\newcommand{\lra}{\longrightarrow}
\usepackage{graphicx}

\title{The Spectrum of Hypersurface Singularities}
\author{Duco van Straten}

\begin{document}
\begin{abstract}
This text is the write-up of a series of lectures 
on the asymptotic mixed Hodge theory of isolated hypersurface 
singularities, held at the Third Latin American school on Algebraic
Geometry and its applications (ELGA 3) in Guanajuato, Mexico, in august 2017.
Its focus is on the classical application of the semi-continuity of the
spectrum due to Varchenko and Steenbrink to the problem of bounding the 
possible singularities on a projective hypersurface.
\end{abstract}
\maketitle

\section{Lecture 1, Monday august 7, 2017}
\subsection{Motivation from classical algebraic geometry}
A hypersurface $Z$ of degree $d$ in projective space $\P^n$ is given as
solution set $\{F=0\}$ of a single homogeneous polynomial 
\[ F(x_0,x_1,\ldots,x_n) \in k[x_0,x_1,\ldots,x_n]_d\] 
of degree $d$. For curves in $\P^2$ and surfaces in $\P^3$ and $k=\R$ one
can make nice pictures of these varieties, and their beauty has inspired 
generations of mathematicians. In the words of {\sc A. Clebsch}:\\

\begin{center} 

{\bf \em ``Es ist die Freude an der Gestalt in einem h\"oheren Sinne, die den Geometer ausmacht.''}
\footnote{{\em ``It is the joy of shape in a higher sense that makes the geometer'', written by A. Clebsch in his obituary for J. Pl\"ucker, {\em Zum Ged\"achtnis an Julius Pl\"ucker}, Dieterichschen Buchhandlung, G\"ottingen, (1872).\\
}}
\vskip 20pt

\includegraphics[height=6cm]{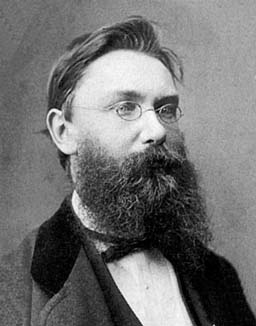}
\hskip 1cm
\includegraphics[height=6cm]{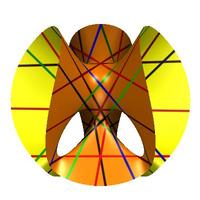}
\vskip 10pt
{\em Alfred Clebsch (1833-1872) and his cubic surface.}\\

\end{center}

For a general choice of coefficients the hypersurface $Z$ is smooth, but
for special choice of coefficients the $n+1$ equations
\[\partial_0F=\partial_1F=\ldots=\partial_n F =0,\;\;\partial_iF:=\partial F/\partial x_i \]
may aquire one or more solutions and define the singular set $\Sigma$ of
$Z$. The question what types of singularities may occur on a hypersurface
of degree $d$ is known only in a very restricted number of cases. For
cubic curves this goes back essentially to {\sc Newton}, for cubic surfaces to
{\sc Schl\"afli}, who listed all possible types of singularities.
The simplest type of singularity is the {\em node}, also called 
{\em ordinary double point} or $A_1$-singularity. 

Such singularities appear as isolated points on the hypersurface and one 
may ask:\\

\centerline{\bf \em Question:}
\vskip 20pt
\begin{center}
{\bf \em What is the maximal number $\mu_n(d)$ of $A_1$-singularities that
can occur on a degree $d$ hypersurface in $\P^n$?}
\end{center}
\vskip 10pt

The value of $\mu_n(d)$ is known only for few values of $n$ and $d$.\\

Let us first look at the case $n=2$, the case of curves of degree $d$ in $\P^2$.
If $C \subset \P^2$ is an irreducible plane curve of degree $d$ and 
\[ n:\widetilde{C} \to C\] 
the normalisation map, then the genus of $\widetilde{C}$ is given by the formula
\[ g(\widetilde{C})= \frac{(d-1)(d-2)}{2}-\delta,\]
where 
\[ \delta=\sum_{p \in C} \delta(C,p).\]
Here $\delta(C,p)$ is the so-called {\em $\delta$-invariant} of the curve
singularity $(C,p)$. Classically it is known under the name of
{\em virtual number of double points}, as it can be shown to be equal
to the number of $A_1$-singularities that emerge from $(C,p)$ if we
perturb the normalisation map slightly.
As $g(\widetilde{C}) \ge 0$, we clearly have
\[\# A_1-\textup{singularities on}\;\; C \le \frac{(d-1)(d-2)}{2} .\]
Equality happens precisely if the normalisation $\widetilde{C}$ is rational.
However, if we allow $C$ to be {\em reducible}, more double points can be
created, and it is not hard to see that the maximal number is realised by
curves $C$ that are union of $d$ lines in general position:
\[ \mu_2(d) =\frac{d(d-1)}{2} .\]

The case $n=3$ of surfaces in three-space is rooted deeply in classical 
algebraic geometry and was a major research direction in italian algebraic 
geometry. Here the problem is incomparably more difficult than for curves 
and the number $\mu_3(d)$ is currently known only for $d \le 6$.

\[
\begin{array}{|c|c|c|c|c|c|c|c|}
\hline
d&1&2&3&4&5&6&7\\
\hline
\mu_3(d)&0&1&4&16&31&65&99 \le \mu_3(d) \le 104\\
\hline
  \end{array}
\]
The is a rich corpus of algebraic geometry around the surfaces with many 
singularities.

\begin{center}
\includegraphics[height=5cm]{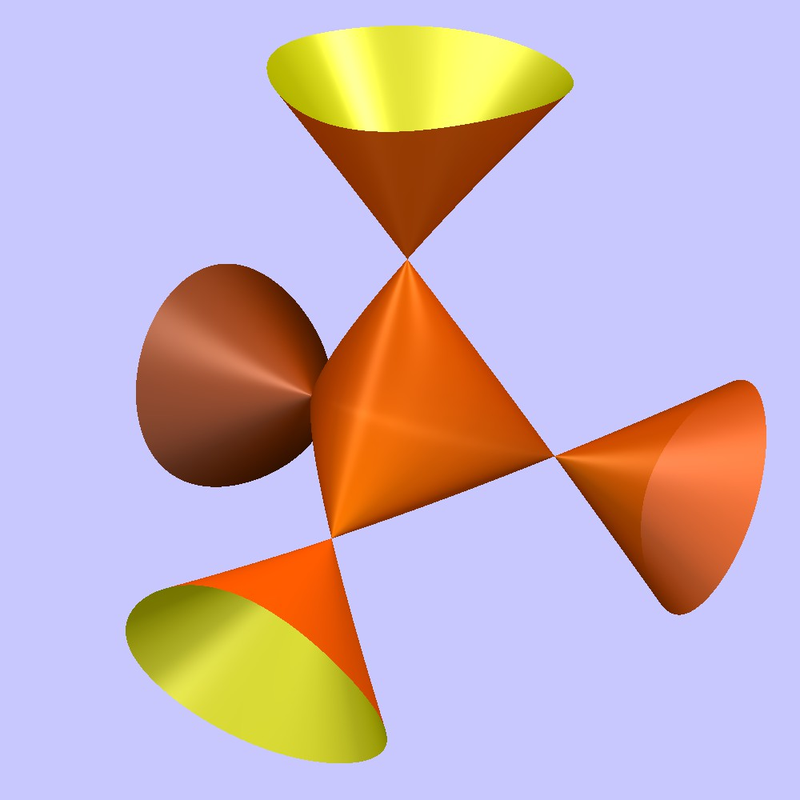}
\hspace{1cm}
\includegraphics[height=5cm]{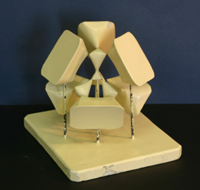}

{\bf \em Cayleys four nodal cubic and Kummers $16$-nodal quartic.}
\end{center}

In $1906$ {\sc A.B. Basset} gave the general upper bound 
\[\mu_3(d) \le \frac{1}{2}(d(d-1)^2-5-\sqrt{d(d-1)(3d-14)+25)}) \]
which was for a long time the best result known.

The bound of {\sc Miyaoka} and {\sc Yau}
\[\mu_3(d) \le \frac{4}{9}d(d-1)^2, \;\;\;(d \ge 4)\] 
for the number of nodes on a surface of general type is the best currently 
known for degree $d \ge 14$.\\ 

\subsection{Spectral upper bound}
In $1981$ {\sc J. B. Bruce} used results on the topology of hypersurface 
singularities to give a new upper bound for the number of nodes on a 
projective hypersurface. Probably inspired by this result,  
{\sc V. I. Arnol'd} formulated in the same year a general upper bound for 
$\mu_n(d)$ in terms of the counting of lattice points in a hypercube between 
two hyperplanes orthogonal to the long diagonal. This conjecture was proven 
by {\sc A. Varchenko}.\\

{\bf \em Theorem ({\sc A. Varchenko}, 1983):}
\[ \mu_n(d) \le A_n(d)\;,\]
where 
\[A_n(d):=\#\left\{(k_1,k_2,\ldots,k_n) \in (0,d)^n \cap \Z^n\;|\;\frac{1}{2}nd-d+1 < \sum_{i=1}^n k_i \le \frac{1}{2} n d \right\} \]
is the so-called {\em Arnol'd number}.
It is elementary to give closed formulas for $A_n(d)$ for fixed $n$.
For example
\[A_3(d)=\left\{ 
\begin{array}{ll} \frac{23}{48}d^3-\frac{9}{8}d^2+\frac{5}{6}d&d \equiv 0 \mod 2\\
\frac{23}{48}d^3-\frac{23}{16}d^2+\frac{73}{48}d-\frac{9}{16}&d \equiv 1 \mod 2,
\\
\end{array}
\right.
\]
so for large $d$
\[ A_3(d) \sim \frac{23}{48}d^3 .\]

{\sc Varchenko}'s bound comes from properties of the {\em  mixed Hodge structure} on the vanishing cohomology of a hypersurface singularity, more specifically from the properties of the {\em spectrum} of such singularities, a notion introduced by
{\sc J. Steenbrink} in $1976$. The definition and proofs involve a mix of topology, algebra and analysis and will be the subject of these lectures.

\subsection{Hypersurface singularities}

We will consider the ring of germs of holomorphic functions
\[ \C\{x_0,x_1,\ldots,x_n\}=\mathcal{O}_{(\C^{n+1},0)}=:\mathcal{O}_{n+1}=:\mathcal{O} .\]
It is a local $\C$-algebra, with maximal ideal 
\[\mathcal{M}=(x_0,x_1,\ldots,x_n) \subset \mathcal{O} .\]
The powers of $\mathcal{M}$ form a descending filtration on $\mathcal{O}$:
\[\mathcal{O} \supset \mathcal{M} \supset \mathcal{M}^2 \supset \ldots \supset \mathcal{M}^k \supset \ldots \supset (0) .\]
The series $f \in \mathcal{M}^k$ are precisely those with vanishing
Taylor series up to order $k-1$.
If $ f \in \mathcal{M}^2$, the series has no linear term and thus has a {\em
critical point} at the origin. We will refer to all series $f \in \mathcal{M}$
as a {\em singularity}.

The automorphisms of the algebra $\mathcal{O}$ can be identified with the
group of coordinate transformations. One says that $f,\;g \in \mathcal{O}$ 
are {\em right equivalent}, notation $f \sim g$, if $f$ and $g$ differ by 
a coordinate transformation. For example
\[x^n+x^{n+1}=(x \sqrt[n]{1+x})^n \sim x^n ,\]
where we use the coordinate transformation $x \mapsto x\sqrt[n]{1+x}$.
If $f \in \mathcal{M}$, but $f \not \in \mathcal{M}^2$, then the implicit
function theorem implies
\[f \sim x_1 .\]
The {\em Morse-lemma} states that if $\phi \in \mathcal{M}^3$, then 
\[x_0^2+x_1^2+\ldots+x_n^2+\phi \sim x_0^2+x_1^2+\ldots+x_n^2 .\]
The classification of functions up to right equivalence starts with the
A-D-E list, and has been pushed to much further. But 
of course the complexity increases without bound, and there is no way 
in which such classification can ever be finished.\\

The {\em Jacobian ideal} of $f \in \mathcal{O}$ is the ideal
\[ J_f=( \partial_0f,\partial_1f,\ldots,\partial_n f) \subset \mathcal{O},\;\; \partial_i f:=\frac{\partial f}{\partial x_i} .\]
The {\em Milnor algebra} of $f$ is the factor ring
\[ \mathcal{O}/J_f\]
and the {\em Milnor number} $\mu(f)$ is defined to be 
\[ \mu(f):=\dim \mathcal{O}/J_f .\]
If \[  \mu(f) < \infty,\]
then one says that $f$ {\em defines an isolated singularity}. As an example, take $f=x^3+y^4$. Its Jacobian ideal is 
$(x^2,y^3)$. So the Milnor algebra has a basis consisting of
(the classes of) the monomials
\[ x^i y^j,\;\;i=0,1,\;\;j=0,1,2\]
so $\mu(f)=6$.
Note that $\mu(f)$ is an {\em invariant} of a singularity: if
$ f \sim g$, then $\mu(f) =\mu(g)$.\\

A fundamental result in singularity theory is the\\

{\em Finite Determinacy Theorem}:\\
{\em If $\mu=\mu(f) < \infty$ and $\phi \in \mathcal{M}^{\mu+2}$, then
\[ f +\phi \sim f .\]
}
In particular, any isolated singularity is right equivalent to a {\em polynomial}! \\
 
A {\em Brieskorn-Pham singularity} is a singularity of the form
\[ f= x_0^{a_0}+x_1^{a_1}+x_2^{a_2}+\ldots+x_n^{a_n},\]
where $a_i \ge 2$. Its Milnor-algebra has a basis of monomials
of the form
\[ x_0^{k_0}x_1^{k_1}\ldots x_n^{k_n}, \;\;0 \le k_i \le a_i-2,\]
and hence
\[ \mu(f) = (a_0-1)(a_1-1)\ldots (a_n-1) .\]
In particular, if all $a_i$ are equal to $d$, the Milnor algebra has
a basis consisting of the monomials inside a $(n+1)$-dimensional cube.
 
\subsection{Spectrum of a hypersurface singularity}

We will give a precise definition of the spectrum later. Here we will give
the main properties that allows one to compute the spectrum of a wide range
of singularities. 

\begin{enumerate}

\item The spectrum $sp(f)$ is a (multi)-set of $\mu:=\mu(f)$ rational numbers
called {\em spectral numbers}:
\[sp(f)=\{\alpha_1 \le \alpha_2 \le \alpha_3 \le \ldots \le \alpha_{\mu}\} .\]
Sometimes it is convenient to pack the information of the
spectrum in a {\em spectral polynomial}:
\[ Sp(f)=\sum_{\alpha \in \Q} n_{\alpha} s^{\alpha} \in \Z[s^{\alpha},\alpha \in \Q],\]
where $n_{\alpha}$ denotes the multiplicity with which $\alpha$ appears in the
spectrum of $f$. Sometimes the spectral numbers are called the {\em exponents}
of a singularity.\\

\item The spectrum is an {\em invariant of a singularity}: if $f \sim g$, then 
$sp(f)=sp(g)$. In fact, something stronger is true: 
if $f \sim u.g$, where $u$ is a unit of $\mathcal{O}$, 
we say that $f$ and $g$ are {\em contact equivalent} and then $sp(f)=sp(g)$.\\

\item {\em Range:} $sp(f) \subset (0,n+1)$.\\
\item {\em Symmetry:} $\alpha_i+\alpha_{\mu-i}= n+1$.\\
\item {\em Thom-Sebastiani:} If $f \in \C\{x_0,x_1,\ldots,x_n\}$ and
$g \in \C\{y_0,y_1,\ldots,y_m\}$ are two series in separate sets of
variables, then 
\[f \oplus g =f(x_0,x_1,\ldots,x_n) +g(y_0,y_1,\ldots,y_m) \in \C\{x_0,\ldots,y_m\}\] 
is called the {\em Thom-Sebastiani sum} of $f$ and $g$.
Then:
\[ sp(f \oplus g) =\{ \alpha+\beta\;|\;\alpha \in sp(f), \beta \in sp(g)\},\]
which can be expressed in term of the spectral polynomials as
\[Sp(f \oplus g) =Sp(f) \cdot Sp(g)\]
\item $sp(x^m)=\{ \frac{1}{m},\frac{2}{m},\frac{3}{m},\ldots,\frac{m-1}{m} .\}$

\end{enumerate}
\vskip 10pt
From these rules one can compute the spectrum for any Brieskorn-Pham singularty.\\

{\em Examples:}\\
\begin{enumerate}
\item $sp(x^2)=\left\{ \frac{1}{2}\right \}$, so 
\[sp(x_0^2+x_1^2+\ldots+x_n^2)=\left\{ \frac{n+1}{2} \right\}.\]
The spectrum consists of a single point at the center of the range $(0,n+1)$.

\item
\[sp(x^2+y^3)=\left\{\frac{1}{2}+\frac{1}{3},\frac{1}{2}+\frac{2}{3}\right\}=\left\{\frac{5}{6},\frac{7}{6}\right\}.\]

\item
\[sp(x^4+y^4)=\left\{ \frac{2}{4},\frac{3}{4},\frac{3}{4}, \frac{4}{4},
\frac{4}{4},\frac{4}{4},
\frac{5}{4},\frac{5}{4},\frac{6}{4} \right\} .\]
Here the spectral numbers appear with a multiplicity, which one
often writes above above the spectral number:

\[
\begin{array}{|ccccc|}
\hline
1&2&3&2&1\\
\hline
\frac{2}{4}&\frac{3}{4}&\frac{4}{4}&\frac{5}{4}&\frac{6}{4}\\[1mm]
\hline
\end{array}
\]

\item
In general, for a Brieskorn-Pham singularity
\[ f=x_0^{a_0}+x_1^{a_1}+\ldots x_n^{a_n},\] 
we can assign to each monomial 
\[ x_0^{k_0} x_1^{k_1} \ldots x_n^{k_n}, \;\;\; 0 \le k_i \le a_i-2,\]
of the Milnor-algebra $\mathcal{O}/J_f$ the number
\[ \sum_{i=0}^n \frac{k_i+1}{a_i} .\]
The (multi)-set of numbers so obtained is precisely $sp(f)$.\\

In fact, as we will see, it is more natural to interpret the number
$\sum_{i=0}^n \frac{k_i+1}{a_i}$ as the quasihomogeneous weight of the
corresponding {\em differential $(n+1)$-form}
\[ x_0^{k_0} x_1^{k_1} \ldots x_n^{k_n}dx_0\wedge dx_1\wedge \ldots \wedge dx_n, \]
where we put 
\[ w(x_i) =\frac{1}{a_i},\]
so that $w(f)=1$.
\end{enumerate}
\vskip 10pt
Recall that more generally a polynomial $f$ is called {\em quasi-homogeneous with weights $w_0, w_1, w_2,\ldots,w_n$}, $w_i \in \Q_{>0}$, if all monomials
$x_0^{\nu_0}x_1^{\nu_1}\ldots x_n^{\nu_n}$ occuring in $f$
are on the hyperplane with the equation 
\[  w_0 \nu_0+w_1\nu_1+\ldots+w_n\nu_n=1 .\]
We call 
\[ w(x_0^{k_0}x_1^{k_1}\ldots x_n^{k_n} dx_0\wedge dx_1 \wedge \ldots \wedge dx_n):=\sum_{i=0}^n w_i(k_i+1)\]
the {\em weight} of the monomial differential form.
We define the {\em Milnor module} as
\[ \Omega_f:=\Omega^{n+1}/df\wedge\Omega^n =\left(\mathcal{O}/J_f\right) \; dx_0\wedge dx_1\wedge \ldots \wedge dx_n .\]
Then the spectral polynomial of $f$ is just the {\em weighted Poincar\'e-series}
of the $\Q$-graded Milnor module 
\[\Omega_f =\bigoplus_{\alpha \in \Q} \Omega_{f,\alpha}\;, \]
that is:
\[Sp(f)=\sum_{\alpha \in \Q} \dim (\Omega_{f,\alpha}) s^{\alpha}.\]

\subsection{Varchenko's theorem} The aim of these lectures is to
describe a proof of the following powerful theorem:\\

{\bf \em Theorem ({\sc Varchenko}, 1983):}\\
{\em Let $Z \subset \P^n$ be a hypersurface of degree $d$ with isolated singular
points $p_1,p_2,\ldots,p_N$. Let the singularity $(Z,p_i)$ be described
by $f_i \in \C\{x_1,x_2,\ldots,x_n\}$. Then for each $\alpha \in \R$ one has
an inequality
\[ \# (\alpha,\alpha+1) \cap sp(x_0^d+x_1^d+\ldots+x_n^d) \ge \sum_{i=1}^N \# (\alpha,\alpha+1) \cap sp(f_i) . \]}
\vskip 8pt
{\em The combined number of spectral numbers for all singularities on a degree
$d$ hypersurface in any open interval of length one is bounded above by the number 
of spectral numbers in the corresponding interval for the singularity $x_0^d+x_1^d+\ldots+x_n^d$.}\\

The spectral numbers of $f= x_0^d+x_1^d+\ldots+x_n^d$ are the weights 
\[ \sum_{i=0}^n \frac{k_i}{d}\]
of the monomials 
\[  x_0^{k_0-1}x_1^{k_1-1}\ldots x_n^{k_n-1} dx_0 \wedge dx_1\wedge \ldots \wedge dx_n,\]
where  
\[  0 < k_i < d .\]
The open interval $(\frac{n-1}{2}+\epsilon,\frac{n+1}{2}+\epsilon)$
contains the spectral number $\frac{n+1}{2}$ of the $A_1$-singularity, so
the number of $A_1$ singularities is bounded by the number of spectral
numbers of $f$ in the corresponding interval, which is exactly the 
Arnold-number $A_n(d)$!\\

{\em Examples:}\\
\begin{enumerate}
\item The Milnor number of $f=x^d+y^d$ is $(d-1)^2$. There are precisely
$d-1$ spectral numbers equal to $1$. So the intervals $(0,1)$ and
$(1,2)$ both contain $((d-1)^2-(d-1))/2=(d-1)(d-2)/2$ spectral numbers.
An interval of length $1$ containing $1$ thus contains at least
$(d-1)+(d-1)(d-2)/2=d(d-1)/2$ spectral numbers. We find that there can be
at most $d(d-1)/2$ nodes on a degree $d$ curve.\\

\item From now we write the multiplicities above the corresponding spectral number:
\[ sp(x^3+y^3+z^3)=\begin{array}{|c|c|c|c|}
\hline
1&3&3&1\\
\hline
 
\frac{3}{3}&\frac{4}{3}&\frac{5}{3}&\frac{6}{3}\\[1mm]
\hline\end{array} .\]
The interval $(\frac{2}{3},\frac{5}{3})$ contains $1+3=4$ spectral numbers, 
so the maximum number of nodes on a cubic surface is at most four. The
Cayley cubic realises this maximum. 

\item
\[ sp(x^4+y^4+z^4)=\begin{array}{|c|c|c|c|c|c|c|}
\hline 
1&3&6&7&6&3&1\\
\hline
\frac{3}{4}&\frac{4}{4}&\frac{5}{4}&\frac{6}{4}&\frac{7}{4}&\frac{8}{4}&\frac{9}{4}\\[1mm]
\hline \end{array}.\]
The interval $(\frac{3}{4},\frac{7}{4})$ contains $3+6+7=16$ spectral numbers, 
so the maximum number of nodes on a quartic surface is at most $16$. Any
Kummer surface realises this maximum.

\item 
\[ sp(x^5+y^5+z^5)=
\begin{array}{|c|c|c|c|c|c|c|c|}
\hline 
1&3&6&10&12&12&10&\ldots\\
\hline 
\frac{3}{5}&\frac{4}{5}&\frac{5}{5}&\frac{6}{5}&\frac{7}{5}&\frac{8}{5}&\frac{9}{5}&\ldots\\[1mm]
\hline 
\end{array} .\]

The interval $(\frac{3}{5},\frac{8}{5})$ contains $3+6+10+12=31$ spectral numbers, so the maximum number of nodes on a quintic surface is at most $31$.
Any Togliatti-surface realises this maximum.\\

\end{enumerate}

{\em Exercises:}\\

1) Work out the spectral bound for the number of $A_2$-singularities on a 
surface of degree $\le 5$.\\

2) Look up the equations for the A-D-E surface singularities and commit them
to memory. Determine their spectra.\\

3) Schl\"afli determined all possible combinations of singularities that can occur on a cubic surface. Check that these combinations are exactly those that
are allowed by the spectral bound. Verify for example the spectrum does not
allow an $A_6$-singularty, but $A_5$ plus an additional $A_1$ is not excluded.\\

3) 
(i) Show that 
\[ x+y+z+u=0,\;\;x^3+y^3+z^3+\frac{1}{4}u^3=0\] 
defines a cubic surface with four nodes ({\em Cayley's cubic}).
Show that 
\[ x+y+z+u+v=0,\;\;x^3+y^3+z^3+u^3+v^3=0\] 
defines a cubic in $\P^4$ with $10$ nodes ({\em Segres cubic}).\\

(ii) Determine the number of nodes of the following cubic in $\P^n$:
\[\sum_{i=0}^{n+1} x_i=0,\;\;\sum_{i=0}^{n+1} x_i^3=0,\;\;\;(n\;\; \textup{even})\]
\[\sum_{i=0}^n x_i=0,\;\;\sum_{i=0}^n x_i^3+\frac{1}{4}x_{n+1}^3=0,\;\;\;(n\;\; \textup{odd})\]

(iii) Show that $A_3(3)=4,\;A_4(3)=10,\;\;A_5(3)=15$ and in general
\[A_n(3) ={n+1 \choose \lfloor \frac{n}{2}\rfloor} \]
and conclude the following result contained in the thesis of {\sc T. Kalker:} 
\[\mu_3(n)=A_3(n) .\] 

4) A general line through a general point $P$ of the Segre cubic $S \subset \P^4$ intersects $S$ in two further points and defines a birational map 
$\phi_P: S \to \P^3$ of degree $2$, that ramifies along a Kummer surfaces, i.e.
a quartic surface with $16$ ordinary double points. Verify these facts.

\section{Lecture 2, tuesday august 8, 2017}

The spectrum has a lot of deeper properties that we will encounter later. 
One of these is the following: if $\alpha \in sp(f)$,
then
\[ \lambda:= \exp(2\pi i \alpha)\]
is an eigenvalue of the cohomological monodromy transformation. In other words,
the spectral numbers $\alpha_1,\alpha_2,\ldots,\alpha_{\mu}$ are {\em specific}
logarithms of the monodromy eigenvalues.\footnote{The spectrum shares this property with the roots of the b-function.} To appreciate this, we first have to
look deeper into the topology of an isolated hypersurface singularity.\\

\subsection{Good representative of a germ}
A power series $f \in \mathcal{M}^2 \subset \mathcal{O}$ with $\mu(f) <\infty$
defines a holomorphic function on a neighbourhood $U$ and by shrinking we may
suppose that $0$ is the only critical point of the function $f:U \to \C$, which
is a singularity of the level set $f^{-1}(0)$.
For $\epsilon$ small enough, the $\epsilon$-ball around $0$  
\[ B_{\epsilon}:=\{ x \in \C^{n+1}\;|\;|x| \le \epsilon\}\]  
is contained in $U$. It is a fundamental non-trivial fact about the analytic 
set $f^{-1}(0)$ that one can find $\epsilon > 0$ with the property that the
boundary $\partial B_{\epsilon'}$ is {\em transverse} to $f^{-1}(0)$ for all
$0 < \epsilon' \le \epsilon$. We now take $0 <\eta$ such that for all 
$t \in S_{\eta}:=\{t \in \C \;|\;|t|\le \eta\} $ the level set $f^{-1}(t)$ is 
transverse to $\partial B_{\epsilon}$. One now puts
\[ X:=X_{\epsilon,\eta}:=B_{\epsilon} \cap f^{-1}(S_{\eta}),\]
\[ S:=S_{\eta} .\]
The function $f:U \to \C$ restricts to a function, again called $f$, 
\[ f: X \lra S\]
\begin{center}
\includegraphics[height=6cm]{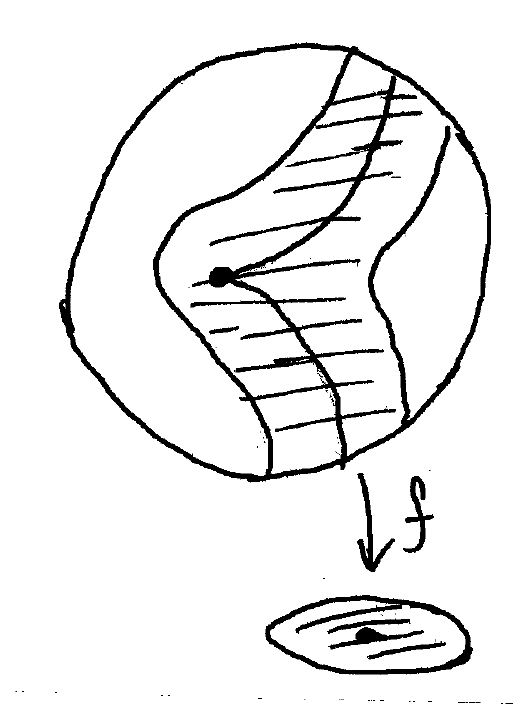}
\end{center}
which is called a {\em good representative of the germ $f \in \mathcal{O}$}.
One puts for $t \in S$:
\[ X_t:=f^{-1}(t) \subset X,\;\;X^*=X \setminus X_0, S^*:=S \setminus\{0\}\]
$f$ restricts to a map
\[ f^*: X^*  \lra S^*\]
that is a locally trivial $C^{\infty}$-fibre bundle by the {\em Ehresmann fibration theorem}. It is called the {\em Milnor-fibration} associated to $f$. 
All fibres $X_t,\; t \in S^*$ are diffeomorphic complex $n$-dimensional 
manifolds with boundary the real $(2n-1)$-dimensional manifold
\[\partial{X_t}=X_t \cap \partial B_{\epsilon} .\]
The manifolds $X_t$ are called the {\em Milnor fibres} of our our 
germ $f \in \mathcal{O}$.
The fibration is trivial near the boundary and consequently all 
$\partial X_t$ can be identified with $L:=\partial X_0$, the {\em link} 
of the singularity. The space $X_0$ is homeomorphic to the topological 
cone over this link, so in particular it is contractible.\\
 
\subsection{Two fundamental theorems} 
For the basic understanding of the topology of hapersurface singularities
the following theorems are of fundamental importance:\\

{\bf \em The bouquet theorem} (Milnor 1968):{\em  Let $\mu:=\mu(f) <\infty$.\\
Then the fibres $X_t$ of the fibration 
$X^* \lra S^*$ have the homotopy type of a bouquet of $\mu$ $n$-spheres:
\[ X_t \equiv \vee_{i=1}^{\mu} S^n.\]}

{\bf \em Corollary:} The homology of the Milnor fibre is non-vanishing
only in degree $0$ and $n$. In fact, for the (reduced) homology in degree 
$n$ one has:
\[\widetilde{H}_n(X_t)=\Z^{\mu} .\]

There is more topological information in the fibre bundle $X^* \lra S^*$. 
As $S^*$ is homotopy equivalent to a circle, the space
$X^*$ can constructed from $X_t \times [0,1]$ by identifying $X_t\times \{0\}$
with $X_t\times \{1\}$ using a glueing map
\[ \tau: X_t \lra X_t .\]
The map $\tau$ can be chosen to be the identity near $\partial X_t$ and is
called the {\em geometric monodromy}. It induces a homological
monodromy transformation 
\[T=\tau_* : H_n(X_t,\Z) \lra H_n(X_t,\Z)\]
which, after a choice of basis, is represented by a $\mu \times \mu$
integral matrix.\\

{\bf \em The monodromy theorem (Brieskorn 1970):} 
The homological monodromy transformation
\[T: H_n(X_t,\Z) \lra H_n(X_t,\Z)\]
is quasi-unipotent. More precisely, there exists a $q \in \N$ such that
\[ (T^q-Id)^{n+1}=0\]
In other words, the eigenvalues of $T$ are roots of unity and the 
Jordan-blocks have size $\le n+1$.\\

Tho monodromy theorem is a general result in algebraic geometry and
various proofs are known. Especially for the case of singularities a
beautiful proof was given by {\sc Brieskorn} in $1970$.\\

{\em Examples:}\\

1) The $A_1$-singularity.\\

If $n=0$ we are dealing with the map $x \mapsto x^2$. The Milnor fibre consists
of two points, which are interchanged under the monodromy. On the reduced
homology group $\widetilde{H}_0(X,\Z)$ the monodromy transformation $T$ is 
multiplication by $-1$.\\
 
If $n=1$ the Milnor fibre $X_t$ defined by $x^2+y^2=t$ has the topology of a 
cylinder; the homology group is generated by the real circle 
\[ \delta(t)=\{(x,y) \in \R^2\;\;|\;\;x^2+y^2=t\} .\]
Projection of the Milnor fibre onto the $x$-line realises it as a two-fold
cover, ramified at the two points $x=\pm\sqrt{t}$
If we run with $t$ once around the origin, these two ramification points
interchange, but the vanishing cycle 'stays the same':
\begin{center}
\includegraphics[height=6cm]{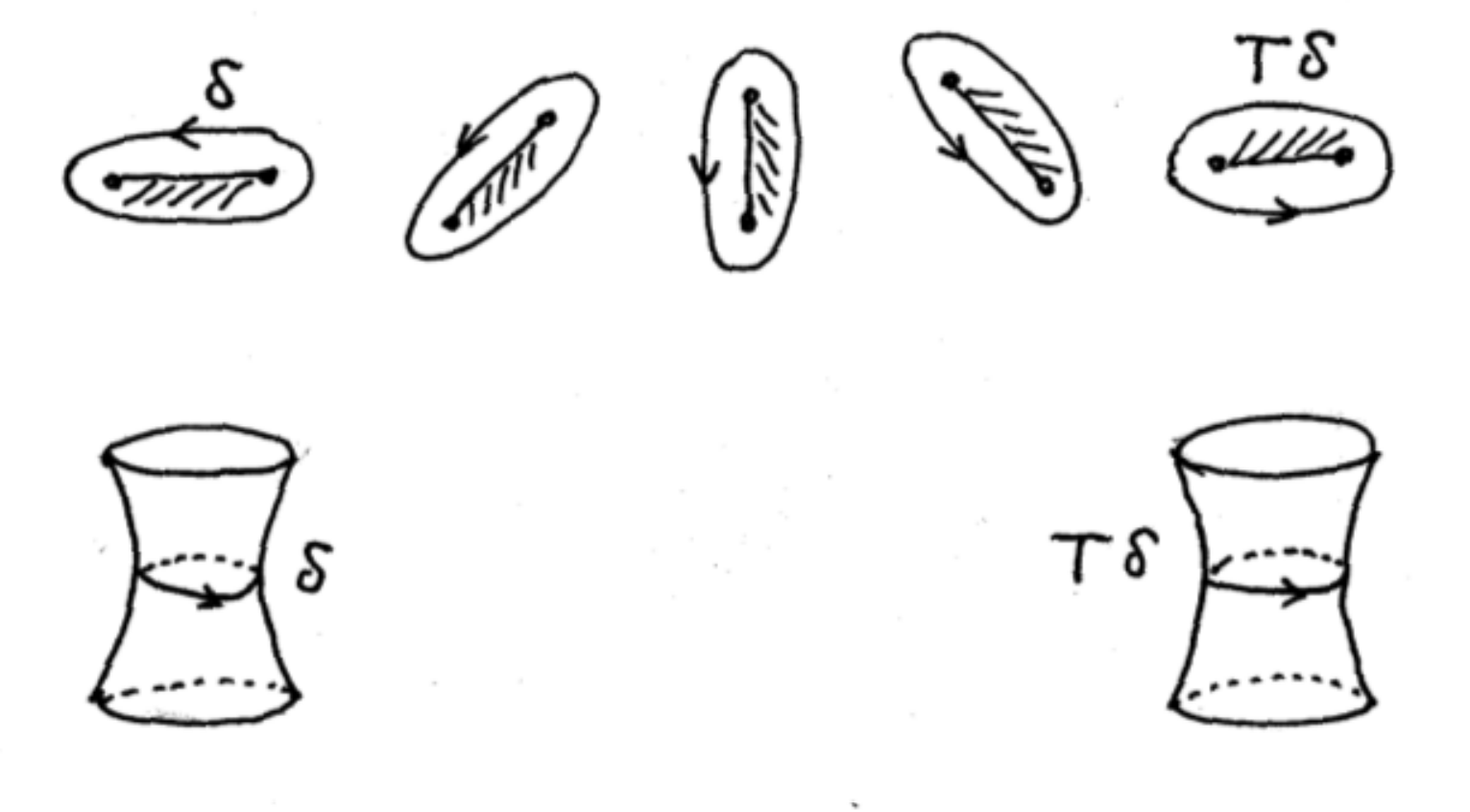}
\end{center}
The monodromy acts as identity on $H_1(X_t,\Z)$.\\

Similarly in higher dimensions: for $n$ even, the monodromy is multiplication
with $-1$, whereas for $n$ odd, it is the identity.\\

However, it is important the realise that even for $n=1$ the monodromy
transformation is geometrically a highly non-trivial map.  
\begin{center}
\includegraphics[height=6cm]{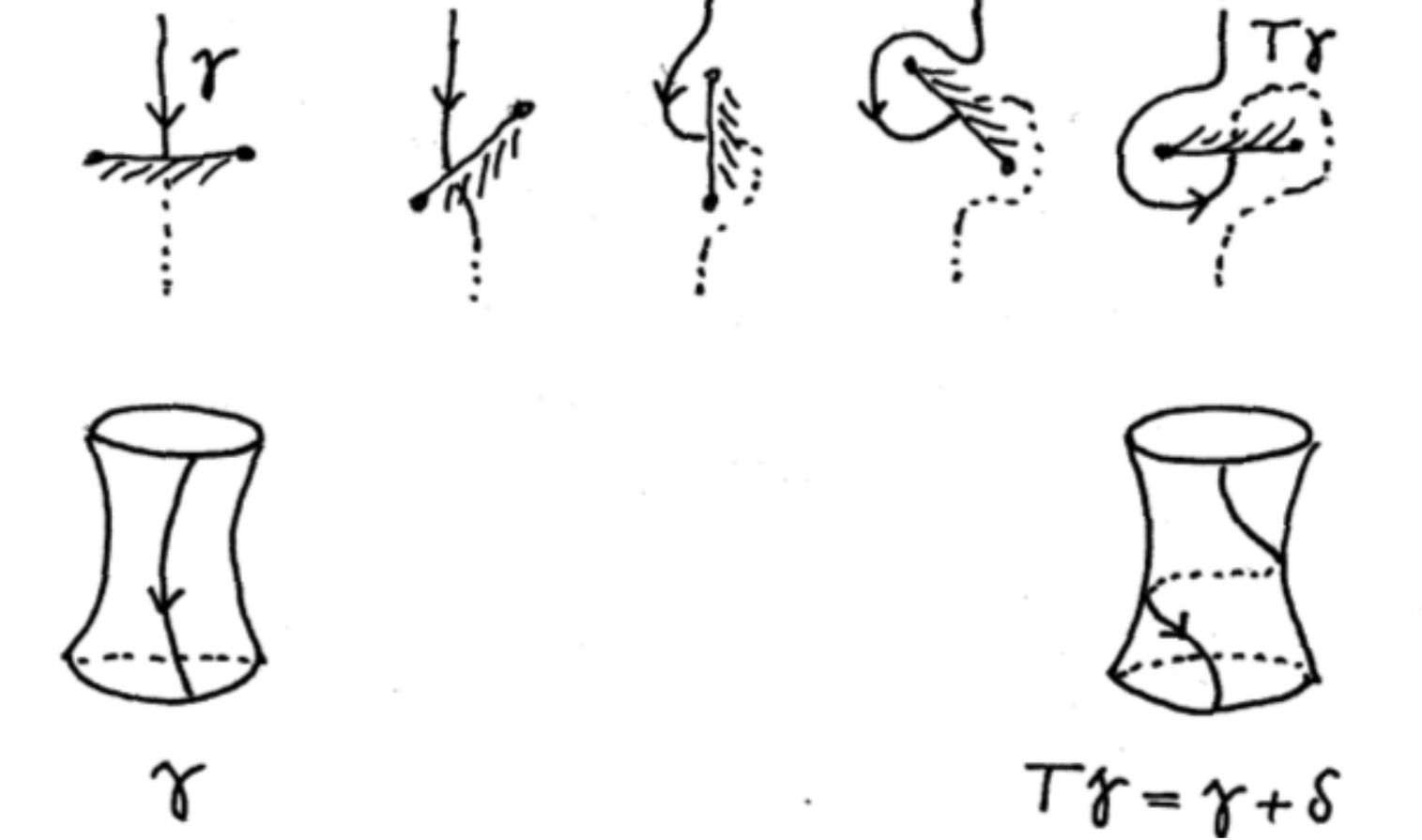}
\end{center}
From this sequence of pictures one sees that the monodromy map $\tau$ is what
in topology is called a {\em Dehn-twist}. Its non-triviality can be detected
by looking at its effect on the relative cycle $\gamma \in H_1(X_t,\partial X_t)$.\\
\footnote{Homologically there is a well-defined map
\[\operatorname{Var}:H_1(X_t,\partial X_t,\Z) \to H_1(X_t,Z),\;\;
[\gamma] \mapsto [T\gamma-\gamma]\]
called the {\em variation mapping}.}

2) $A_2$-singularity $f=y^2+x^3$.
The Milnor fibre $$X_t=\{y^2+x^3=t\}$$ 
has the topology of an elliptic curve, with a disc removed, so indeed has 
the homotopy type of the wedge of two circles. 
To 'see' the monodromy, we perturb the function $f$ and consider
\[f_{\lambda}=y^2+x^3-\lambda x\]
where $\lambda$ is a small positive real number. For a fixed value of 
$\lambda$,  the function $f_{\lambda}$ has two critical points. These points 
appear when the cubic polynomial 
\[ x^3-\lambda x-t=0\] 
has a repeated root. These points are of 
form $(x_{+},0)$ and $x_{-},0)$ with critical values 
\[t_+=\sqrt{\frac{4 \lambda^3}{27}},\;\;\;t_{-}=-\sqrt{\frac{4 \lambda^3}{27}} .\]
We choose  a basis of the homology $\delta_+, \delta_-$ such that if we move 
towards $t_{\pm}$, the cycle $\delta_{\pm}$ is vanishing.
\begin{center}
\includegraphics[height=5cm]{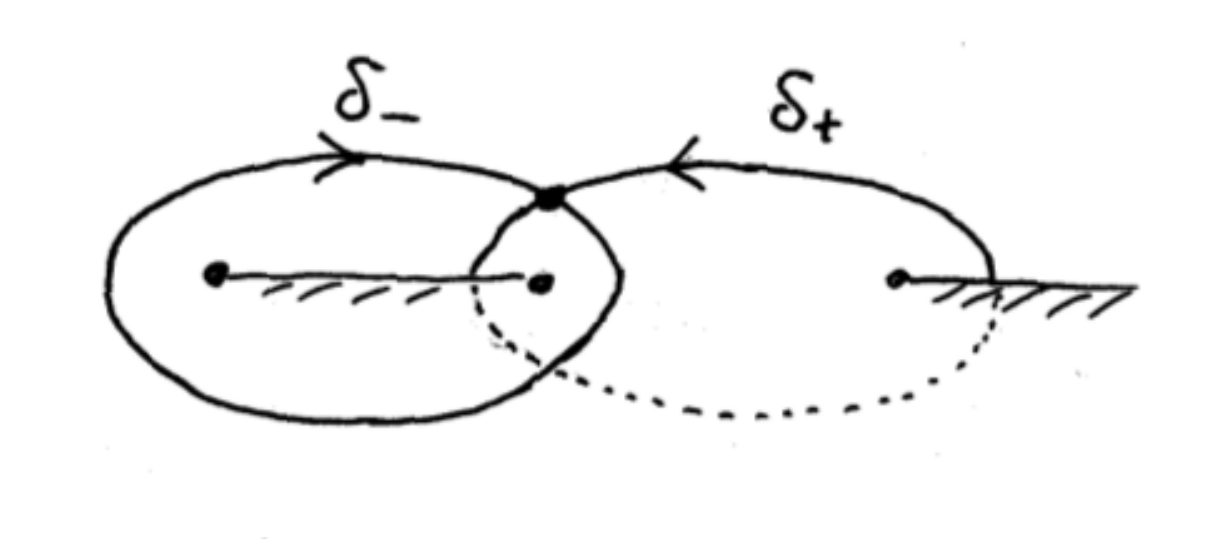}
\end{center}
We orient the vanishing cycles such that the intersection product is given by
\[ <\delta_+,\delta_- > =1,\] 
as indicated in the picture.
We consider in the disc $S_{\lambda} \subset \C$ the straight line paths $\gamma_{\pm}$ from a base point $*$ to $t_{\pm}$. If we start from $*$, follow the path $\gamma_{\pm}$ to a point close to $t_{\pm}$ and encircle it in the positve direction, and then return to the base point $*$ we obtain two local monodromy
transformations $T_{\pm}$. 
\begin{center}
\includegraphics[height=6.5cm]{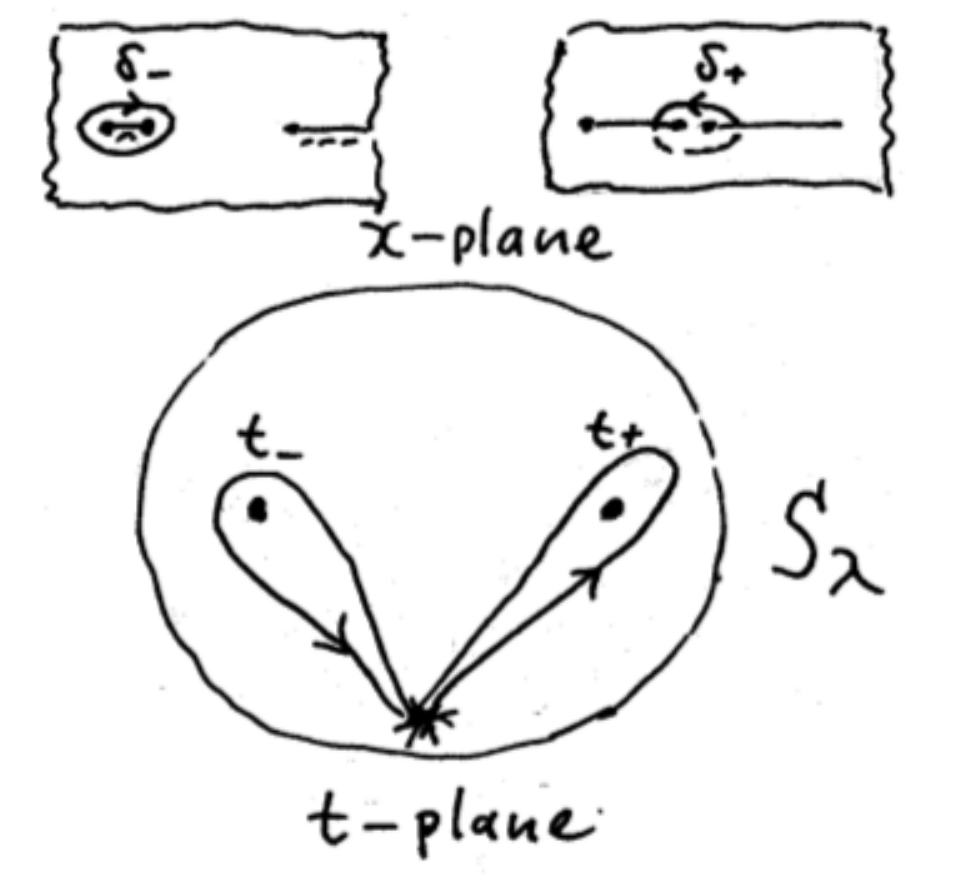}
\end{center}

From the analysis of the $A_1$-case we can
see that
\[T_{+}(\delta_+)=\delta_{+},\;\;T_{+}(\delta_-)=\delta_{-}+\delta_+,\;\;\;T_{-}(\delta_+)=\delta_{+}-\delta_-,\;\;T_{-}(\delta_-)=\delta_{-}\]
Hence, with respect to the basis $\delta_+,\delta_-$ these
transformations $T_+$ and $T_-$ are represented by the following matrices:
\[T_+=\left( \begin{array}{rc}1&1\\0&1\end{array}\right),\;\;\;T_-=\left( \begin{array}{rc}1&0\\-1&1\end{array}\right) .\]
 
The monodromy transformation $T$ of the $A_2$-singularity is obtained by
encircling {\em both} critical points, so is given by the matrix product
\[ T=T_{-} T_{+} =\left( \begin{array}{rc}1&0\\-1&1\end{array}\right) \left( \begin{array}{rc}1&1\\0&1\end{array}\right)=\left( \begin{array}{rc}1&1\\-1&0\end{array}\right)\] 
It is easy to check that 
\[ T^6=I\]
and the eigenvalues of $T$ are the primitive sixth roots of unity
\[ e^{2\pi i/6}\;\; \textup{and}\;\; e^{- 2\pi i/6} .\]

The above idea of using a small perturbation to analyse the topology
of a singularity is useful in many other situations as well. 
Using these ideas that can be traced back to {\sc Lefschetz}, 
{\sc Brieskorn} gave an alternative proof of the bouquet theorem, 
that we will sketch now.\\

\begin{enumerate}
\item We consider a perturbation of our function $f$ of the form
\[ f_{\lambda} =f+\sum_{i=0}^n \lambda_i x_i .\]
If we take the $\lambda_i$ small enough, the function $f_{\lambda}$ 
will have a finite number of critical points $p_1,p_2,\ldots,p_N$
in the set $X_{\lambda}:=f_{\lambda}^{-1}(S_{\eta})$.
As the critical space, defined by the vanishing of the partial derivatives 
$$\partial_0 f_{\lambda}=\partial_1 f_{\lambda}=\ldots=\partial_n f_{\lambda}$$ 
is a complete intersection, it is Cohen-Macaulay, hence the sum of the Milnor 
numbers at the points $p_i,i=1,\ldots,N$ add 
up to the value at $\mu=\mu(f)$ for $\lambda=0$.\\
\item By choosing the $\lambda_i$ not only small, but also {\em generic}, 
the critical points of $f_{\lambda}$ will all be non-degenerate and have different critical values under the map $f_{\lambda}$. In particular, there will be exactly 
$\mu$ critical points and $\mu$ critical values in the disc $S_{\eta}$. 
We obtain what is called a {\em morsification of $f$}.\\
\item Now choose paths $\gamma_{i}, i=1,2,\ldots,\mu$ connecting a base
point $* \in \partial S_{\eta}$ with the critical values $t_1,t_2,\ldots,t_{\mu}$.\\ 

\begin{center}
\includegraphics[height=6cm]{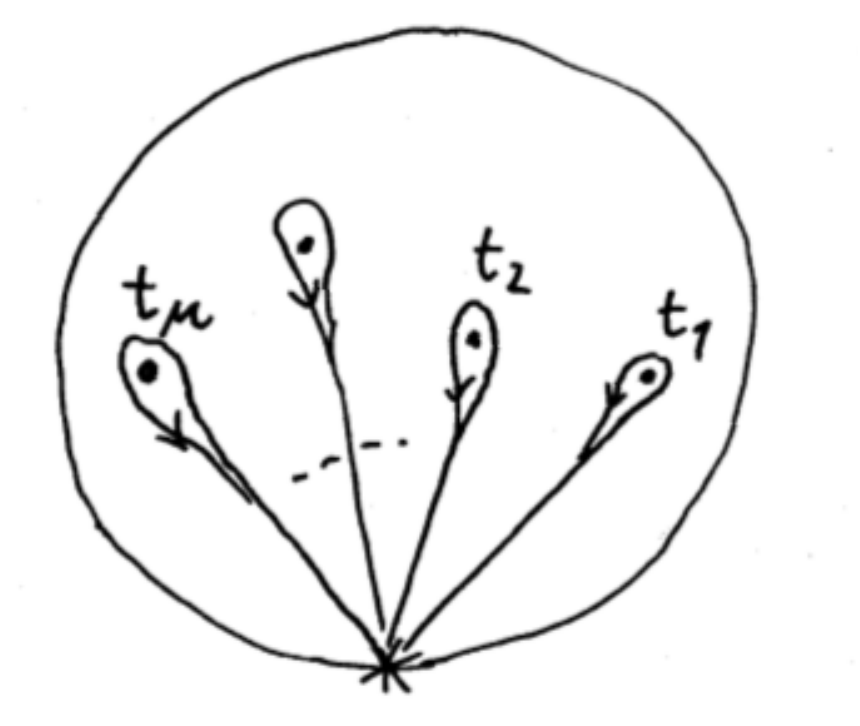}
\end{center}

\item The preimage of the disc $S_{\eta}$ is $X_{\lambda}$, which is contractible.
Furthermore, the disc $S_{\eta}$ contracts to the union of the paths $\gamma_i$, and their preimage under $f_{\lambda}$ contracts to the the fibre $X_{*}$ over
the base point, to which we glue the Lefschetz thimbles $\Gamma_i$, $i=1,2,\ldots,\mu$, lying over the paths $\gamma_i$, $i=1,2,\ldots,\mu$. So we see that
by glueing $\mu$-thimbles to the fibre $X_{*}$ we obtain a contractible space.\\

\begin{center}
\includegraphics[height=6cm]{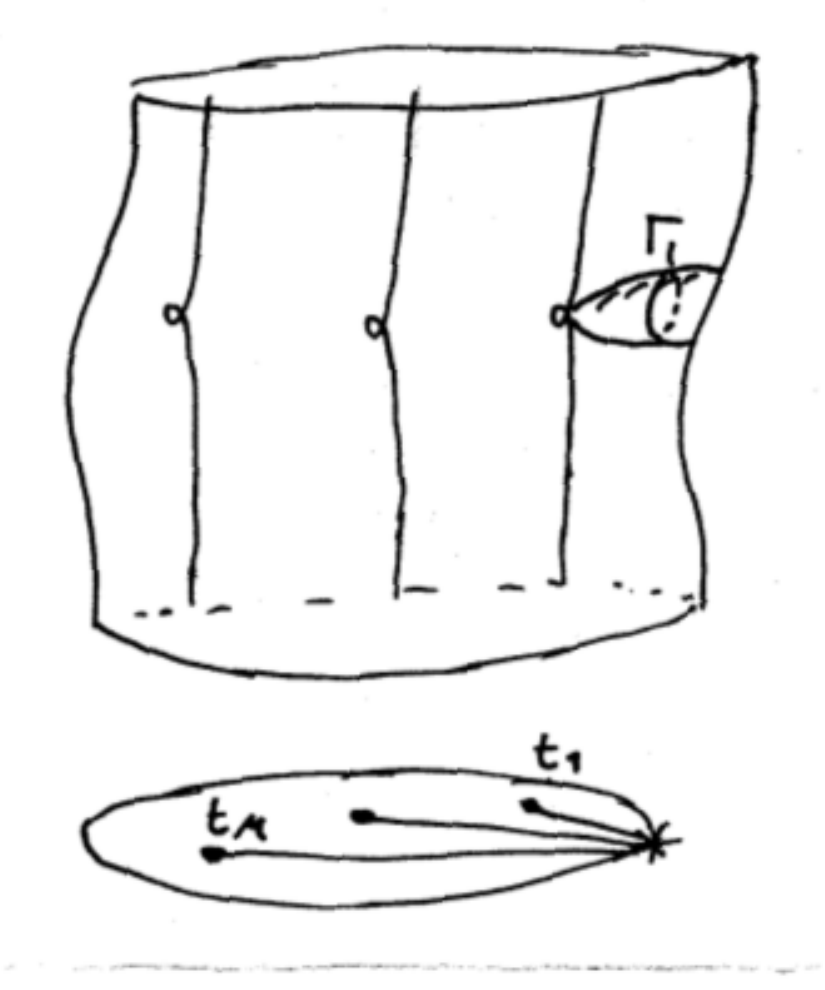}
\end{center}

\item The thimbles are topological discs of real dimension $n$. So by glueing
$\mu$ $n$-discs to $X_{*}$ we obtain a contractible space. If follows from CW-topology that $X_{*}$ can only have the homotopy type of a bouquet of 
$n$-spheres!\\
\end{enumerate}

It is far more difficult to understand the nature of the monodromy transformation from the above morsification picture.\\

{\bf \em Theorem ({\sc F. Pham, 1965}):}\\

{\em  For the singularity
$f=x_0^{a_0}+x_1^{a_1}+\ldots x_n^{a_n}$
with Milnor number 
\[\mu=(a_0-1)(a_1-1)\ldots(a_n-1),\]
the monodromy $$T:H_n(X_t) \lra H_n(X_t)$$ is of finite order
$$e:=lcm(a_0,a_1,\ldots,a_n) .$$
The eigenvalues of $T$ are the numbers
\[ \omega_0 \omega_1 \ldots \omega_n,\]
where the $\omega_i$ are (non-trivial) $a_i$-th roots of unity.}\\ 

Because $f$ is obtained as Thom-Sebastiani sum from functions
depending on a single variable, it is easy to prove the result,
using the fact for a Thom-Sebastiani sum the homology and
monodromy transformation can be identified with the tensor product
of the factors:
\[\widetilde{H}_n(X(f\oplus g))=\widetilde{H}_n(X(f))\otimes \widetilde{H}_n(X(g)),\;\;\; T({f\oplus g})=T(f) \otimes T(g)\]
Furthermore, L\^{e} had shown that the monodromy to be of finite order for any
irreducible curve singularity and the question arose if the monodromy of 
any isolated singularity was always of finite order.

\subsection{Examples of A'Campo and Malgrange}

The answer is: NO ! The first examples were given by {\sc N. A'Campo}
in $1973$. Consider the function
\[ f=(x^2+y^3)(x^3+y^2)=x^2y^2+x^5+y^5+x^3y^3 \sim x^2y^2+x^5+y^5 \]
Its zero-level consist of two transverse cuspidal branches. 
\begin{center}
\includegraphics[height=5cm]{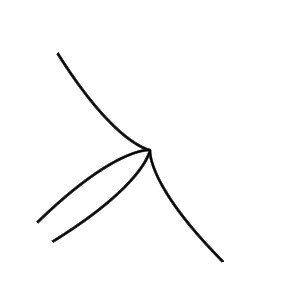}
\end{center}

Let us consider the {\em embedded resolution} of this function. Blowing up the origin introduces
a component, along which the pull-back of $f$ vanishes with order $4=mult_0(f)$.
The strict transform consists of smooth branches, which are tangent to this
divisor. Blowing up at these two points introduces two further exceptional 
curves, along which the function (always pull-back) vanishes with order $4+1=5$.
The strict transforms now are transverse to the first divisor and the newly introduced divisors. Blowing once more at these two point introduces two further curves, along which the pull-back of $f$ vanish with order $1+5+4=10$.
The configuration of exceptional curves now looks as follows. 
\begin{center}
\includegraphics[height=5cm]{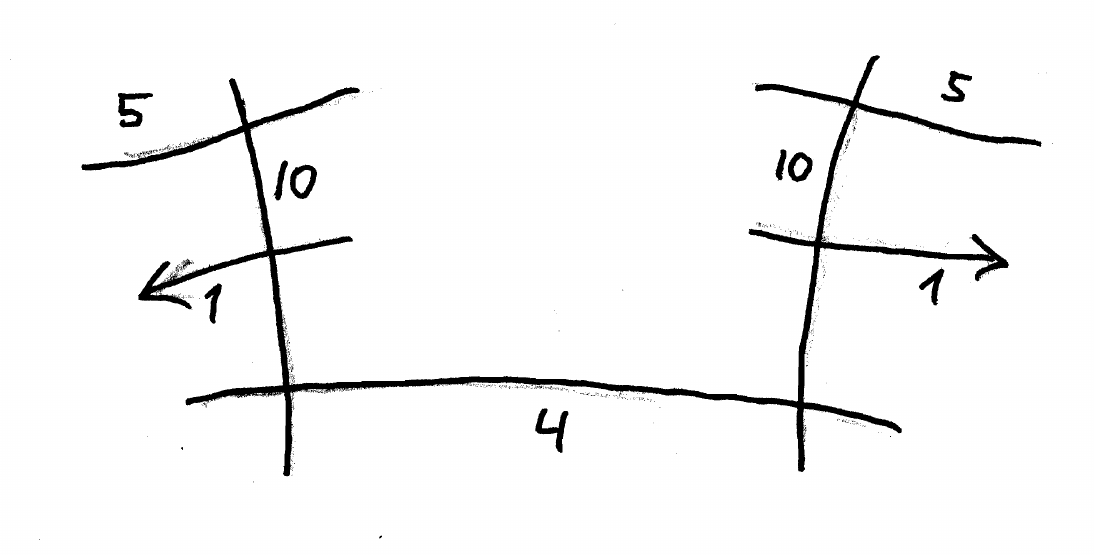}
\end{center}

The strict transforms of the two branches are indicated by an arrow. This
configuration now leads to a very precise model of the Milnor-fibre, obtained
from cyclic coverings of these curves determined by the multiplicities of $f$.
As explained in detail in the book of {\sc Brieskorn} and {\sc Kn\"orrer}, the
result looks as follows:

\begin{center}
\includegraphics[height=5cm]{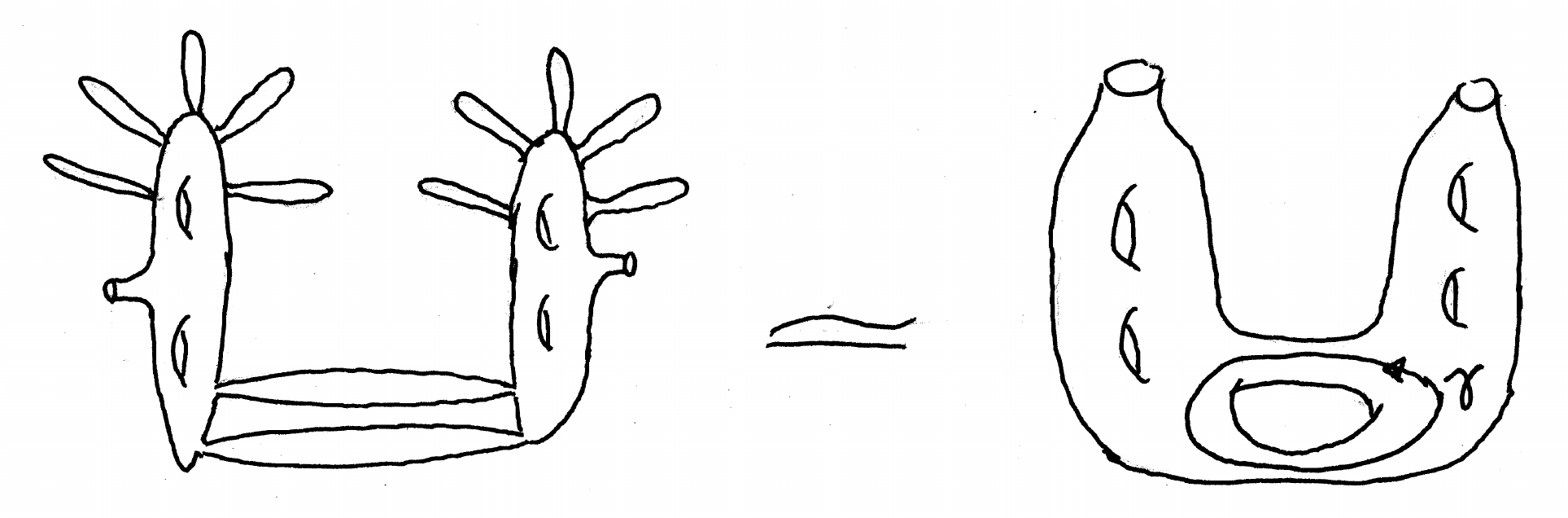}
\end{center}

In this way A'Campo showed that the curve $\gamma$ in the picture satisfies
$(T^{10}-I) \gamma \neq 0$, so one has 
\[(T^{10}-I) \neq 0,\;\;\;(T^{10}-I)^2=0\]
 
In his paper,  {\sc A'Campo} asked if there were examples of isolated singularities in $n+1$ variables for which the monodromy has a Jordan block of size $n+1$. It is hard to analyse higher dimensional examples along these geometrical lines and the question was answered by {\sc Malgrange}, 
using a totally different kind of argument.\\

{\bf Theorem ({\sc B. Malgrange}, 1973):}\\

{\em The monodromy transformation $T$ of the
singularity
\[f=(x_0 x_1\ldots x_n)^2+x_0^{2n+4}+x_1^{2n+4}+\ldots+x_n^{2n+4}\]
has a Jordan block of size $n+1$.}\\

The argument given by {\sc Malgrange} runs as follows. 
For $t$ real and positive, the set
\[ \delta(t):=\{ x \in \R^{n+1}\;|\;f(x)=t\} \subset X_t\]
has the topology of an $n$-dimensional sphere. It is 
the boundary of the set
\[E(t) :=\{ x \in \R^{n+1}\;|\;f(x) \le t\} .\]
Clearly, $E(t)$ is topologically an $(n+1)$-ball; it can be
seen as the Lefschetz thimble of the vanishing cycle $\delta(t)$.
Now consider
\[I(t):=\int_{\delta(t)} x_0dx_1dx_2 \ldots dx_n =\int_{E(t)} dx_0dx_1\ldots dx_n=Vol(E(t)) .\]
{\sc Malgrange} claims that for $t \to 0$, the integral $I(t)$ behaves like
\[ I(t) \sim C t^{1/2} (\log t)^n,\]
where $C$ is a constant.
If we analytically continue $I(t)$ along a path going once around the origin
in the positive direction, $I(t)$ changes its value to
\[ (-t^{1/2})(\log t+2\pi i)^n,\]
and will not return to its original value by going around any number of times.
On the other hand, this is equal to the integral of $\omega=x_0dx_1dx_2\ldots dx_n$ over the transformed cycle $T\delta(t)$. Clearly, $T$ can not be of finite order, and a little  further thinking shows that the $(\log t)^n$-term creates a Jordan block of size $n+1$.\\

How to find the asymptotic expansion of the integral $I(t)$? For this
Malgrange suggestes the following method.
Let 
\[F(t):=\{x \in \R^{n+1}\;|\;\; (x_0x_1\ldots x_n)^2 \le t, x_i^{2n+4}\le t\}.\]
Then it is easy to see that
\[F \left(\frac{t}{n+1}\right) \subset E(f) \subset F(t),\]
so that it is sufficient to study the asymptotic expansion of $Vol(F(t))$,
which is elementary.
To illustrate the point, let us look at the case $n=2$, so we look at
the curve $\delta(t): (xy)^2+x^6+y^6=t$ for various values of t
\begin{center}
\includegraphics[height=4cm]{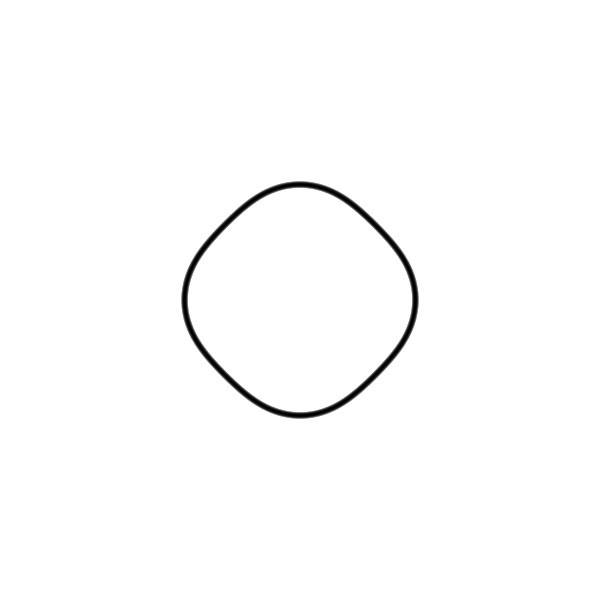}\vspace{1cm}
\includegraphics[height=4cm]{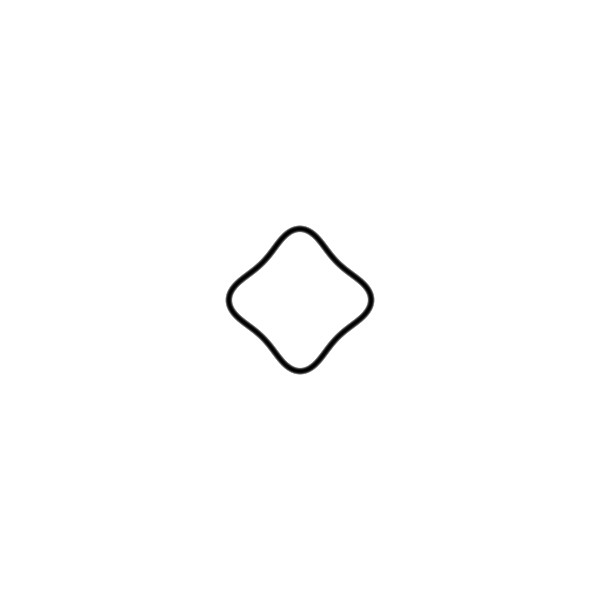}\vspace{1cm}
\includegraphics[height=4cm]{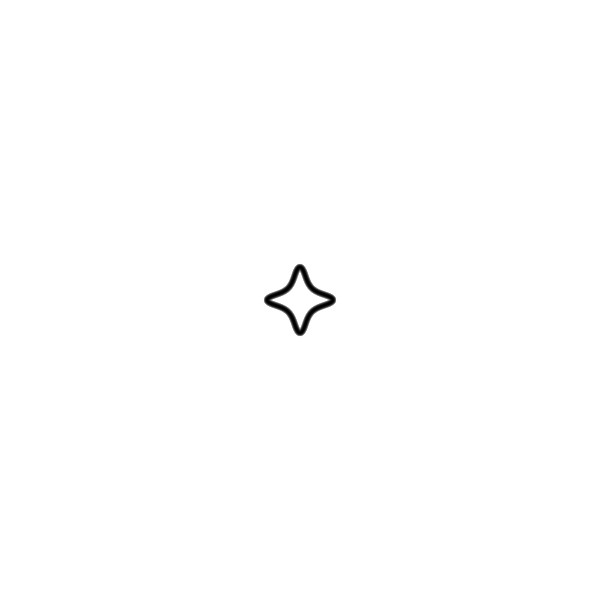}

{\bf \em The vanishing cycle $\delta(t)$ for smaller and smaller values of $t$.}
\end{center}
The volume of the region $F(t)$ is, up to a factor four, equal to the area
of the subset of the quadrant $x \ge 0, y\ge 0$ given by the
inequalities
\[ xy \le t^{1/2},\;\;x\le t^{1/6},\;\;\;y \le t^{1/6}\]
Doing the integral, we  find
\[
\begin{array}{rcl}
Area(F(t)/4)&=& t^{1/3}t^{1/6}+\int_{t^{1/3}}1{t^{1/6}} t^{1/2}\frac{dx}{x}\\[2mm]
            &=& t^{1/2}-\frac{1}{6}t^{1/2}\log t .\\
\end{array}
\]
So just by looking at the shape of $\delta(t)$ for small $t$ is is 'clear'
that a logarithm has to appear in the expansion, and hence that the monodromy
is of infinite order!\\

{\em Exercise:}
1) Show that the Milnor fibre of the $A_1$-singularity in dimension $n$
can be identified with the disc bundle inside the cotangent bundle to the
$n$-sphere.\\

2) Give an analysis of the monodromy of the  $A_3$-singularity using a
morsification. Determine the three local monodromies and compute their product.
What is the order of the resulting monodromy transformation?\\

3)  Give an analysis of the Milnor-fibre of the $A_2$-singularity
using an embedded resolution. Show that one obtaines the elliptic
curve with an automorphism of order $6$.\\ 

4) Work out the details for the expansion of the integrals considered
by {\sc Malgrange}. Can you find the sub-leading order?

\section{Lecture 3, wednesday, august 9}
We saw yesterday in an example that it was useful to consider
integrals of differential forms over vanishing cycles. 
Such integrals turn out to expand in a series of the form
\[ I(t) =\int_{\delta(t)}\omega=C t^{\alpha} \log t^k+\ldots\]
Such an integral represents a multi-valued function on the punctured disc
$S^*$. If we fix a determination of this multi-valued function
near the point $t_0 \in \partial S$
and continue along the path $t_0 e^{i\phi}$, where the argument
$\phi$ runs from $0$ to $2\pi$, the function $I(t)$ changes to
$$T(I(t))=I(t\exp(2\pi i)) .$$
Note that we will have the relation
\[ T(I(t))=\int_{T\delta(t)} \omega .\] 
It is easy to analyse what happens with the terms in the expansion:
upon analytic continuation, 
the function $t^{\alpha}$ changes to $t^{\alpha+2\pi i \alpha}$, so 
\[T(t^{\alpha}) = e^{2\pi i \alpha} t^{\alpha} .\] 
In words, it gets multiplied by the number $e^{2\pi i \alpha}$, which 
is a root of  unity precisely when $\alpha \in \Q$.
 
The logarithm $\log t$ changes to $\log t +2\pi i$:
$$ T(\log t)= \log t+2\pi i .$$
So we see that the pair of functions $\log t, 1$ create a Jordan-block
of size $2$:
\[
T \left(\begin{array}{c}\log t\\ 1 \end{array} \right)=\left(
\begin{array}{cc}
1&2\pi i\\
0&1\\
\end{array}
 \right)\left(\begin{array}{c} \log t \\ 1 \end{array} \right)\]
and
\[(T-1)(\log t) \neq 0,\;\;\;(T-1)^2(\log t)=0 .\] 
Similarly, 
\[
T \left(\begin{array}{c}(\log t)^2\\\log t\\1 \end{array} \right)=\left(
\begin{array}{ccc}
1&4\pi i&(2\pi i)^2\\
0&1&4\pi i\\
0&0&1\\
\end{array}
 \right)\left(\begin{array}{c} (\log t)^2 \\\log t\\ 1 \end{array} \right)\]
so,
\[ (T-1)^2(\log t)^2 \neq 0,\;\;\;(T-1)^3(\log t)^2=0\]
etcetera.

In general, if $\alpha=k/q$, then one has

\[(T^q-I)^{n+1} t^{\alpha} (\log t)^n =0.\]
 
\subsection{What integrals to consider?}
First, let us recall  some aspects of the Poincar\'e residue. If
$Z$ is a complex manifold and $V$ a smooth divisor defined by a 
holomorphic function $f$, we say that a holomorphic differential 
form $\eta$ is {\em residue} of a meromorphic form $\omega$
on $Z$, if we can write
\[ \omega=\eta \wedge \frac{df}{f} + \omega',\]
where $\omega'$ is regular on $Z$. One writes
\[Res_V(\omega)=\eta_{|V}.\] 
The residue theorem from the elementary theory of holomorphic functions
has the following generalisation:
\[ \int_{\delta} Res_V(\omega) =\frac{1}{2\pi i}\int_{\Delta} \omega\]
Here $\Delta \subset Z\setminus V$ is obtained from $\delta$ by applying
the tube operation: we replace each point of the cycle $\delta$ by
a small circle encircling the divisor $V$ in the positive direction.\\

We will apply this in particular to the situation that arises in the
geometry of a good representative of an isolated hypersurface singularity.
So, we consider a holomorphic  $(n+1)$-form on $X$:
\[  \omega \in \Gamma(X,\Omega_X^{n+1}) .\]
Let us fix $t \in S^*$ and consider the following differential form
on the Milnor fibre $X_t$:
\[ \eta_t :=Res\left(\frac{\omega}{f-t}\right) \in \Gamma(X_t,\Omega_{X_t}^n) .\]
The period integrals we want to consider are of the form
\[ I(t):=\int_{\delta(t)} Res\left(\frac{\omega}{f-t}\right) =\frac{1}{2\pi i}\int_{\Delta(t)}
\frac{\omega}{f-t},\]
where $\delta(t) \in H_n(X_t,\Z)$. If we fix a sector $W$ in the punctured disc $S^*$ and use 
parallel transport of a cycle $\delta(t_0) \in H_n(X_{t_0},\Z)$
to obtain cycles $\delta(t) \in H_n(X_t,\Z)$ for $t \in W$, we obtain a well-defined holomorphic function $I(t)$ in $W$.\\

{\em Example:} Let us consider a quasi-homogeneous singularity $f$ with
weights $w_i=w(x_i)$. Put
\[ \omega=x_0^{\nu_0}x_1^{\nu_1}\ldots x_n^{\nu_n} dx_0dx_1\ldots dx_n .\]
For any choice of cycle $\delta(t)$, the integral is seen to have a
scaling property: if we replace $x_i$ be $\lambda^{w_i}x_i$ in the
integral, we find
\[I(t)=\frac{1}{2\pi i}\int_{\Delta(t)} \frac{\lambda^{w(\omega)} \omega}{\lambda f(x)-t}=\lambda^{w(\omega)-1}I(t/\lambda)\]
So putting $t=\lambda$ we get the scaling relation:
\[I(t) =t^{w(\omega)-1} I(1) \]

These period integrals are rather trivial. They have finite monodromy, but
the weight of the differential form is reflected in the speed of vanishing
of the period integral.\\

{\bf \em The regularity theorem:}\\

{\em All period integrals of the form
\[ I(t):=\int_{\delta(t)} Res\left(\frac{\omega}{f-t}\right)\]
attached to an isolated hypersurface singularity have an expansion 
(convergent on any sector $W\subset S^*)$ of the form
\[ I(t)=\sum_{\alpha, k} C_{\alpha,k} t^{\alpha}(\log t)^k\]

Furthermore:

\begin{itemize}
\item $e^{2\pi i \alpha}$ is an eigenvalue of the monodromy $T:H_n(X_t,\Z) \lra H_n(X_t,\Z)$.
\item $0 \le k \le n+1$.
\item $\alpha > -1$.
\end{itemize}
}

The last inequality is of particular importance, as was first
recognised by {\sc Malgrange}. It can be deduced from the following\\

{\bf \em Positivity Lemma ({\sc Malgrange (1974)}):}\\
If $\eta \in \Gamma(X,\Omega^n_X)$ is a global $n$-form on $X$
and if we put 
\[\eta_t :=\eta_{| X_t} \in \Gamma(X_t,\Omega_{X_t})\]
then, if  $t \lra 0$:
\[\int_{\delta(t)} \eta_t \lra 0 .\]

The idea is that the integral of a holomorphic differential 
form over cycles that shrinks to $0$ has to go to zero, just because we
integrating over smaller and smaller sets.\\

\subsection{Connection on the cohomology bundle}
We have seen that the restriction of $f:X \lra S$ to $X^*$ defines a
is a locally trivial fibre bundle over $S^*$. The direct image
\[ R^nf_*\Z_{X^*} \]
is a {\em local system} with fibre over $t \in S^*$ equal to $H^n(X_t,\Z)$.
The monodromy transformation $T$ describes the non-trivial twisting in this
lattice bundle. We will denote its complexification by $H_S^*$:
\[ H_{S^*}=\C \otimes_{\Z} R^nf_*\Z_{X^*}=R^nf_*\C_{X^*}\]
One may associate to this local system of $\mu$-dimensional vector spaces 
a vector bundle of rank $\mu$ by tensoring this with the holomorphic 
functions on $S^*$:
\[ \mathcal{H}_{S^*}:=\mathcal{O}_S \otimes_{\C} H_{S^*}\] 
which we will call the {\em cohomology bundle}. Its sections describe
cohomology classes that depend holomorphically on $t \in S^*$.
On the vector bundle $\mathcal{H}_{S^*}$ there is a tautological connnection
\[ \nabla: \mathcal{H}_{S^*} \lra \Omega^1_{S^*} \otimes \mathcal{H}_{S^*},\]
defined by
\[ \nabla( g \otimes h) = dg \otimes h .\]
If we denote by $t$ the local coordinate on $S$, then we have the
standard vector field $\frac{\partial}{\partial t}$ dual to $dt$ and we get
a map
\[ \partial_t:= \nabla_{\partial/\partial_t}:\;\mathcal{H}_{S^*} \lra \mathcal{H}_{S^*},
\;\;\;g \otimes h \mapsto \frac{dg}{dt}\otimes h .
\]
The local system of {\em horizontal sections} of the 
cohomological bundle is identified with the local system we started with:
\[ H_{S^*}=Ker(\partial_t) \subset \mathcal{H}_{S^*} .\]

One would like to extend $\mathcal{H}_{S^*}$ to a vector bundle $\mathcal{H}_S$
on $S \supset S^{*}$. This is easy and can be done in many possible ways. 
As $S^*$ is a Stein space, the vector bundle $\mathcal{H}_{S^*}$ is in fact 
a {\em trivial} vector bundle. We can pick {\em any} basis of sections 
$e_1,e_2,\ldots,e_{\mu} \in \Gamma(S^*,\mathcal{H}_{S^*})$ and consider the
vector bundle
\[ \mathcal{H}_S:=\oplus_{i=1}^{\mu} \mathcal{O}_S e_i \subset j_* \mathcal{H}_{S^*},\]
where $j: S^* \lra S$ is the inclusion map.\\

 What happens to the connection? Well, we can write out the effect of
$\partial_t$ on the basis sections $e_i$:
\[ \partial_t e_j =\sum_{i=1}^{\mu} A_{ij}(t) e_i,\]
and obtain a matrix of holomorphic function on $S^*$
\[ A(t)=(A_{ij}(t)) \in \operatorname{Mat}(\mu \times \mu,\mathcal{O}_{S^*}) ,\]
called the connection matrix. A section 
\[ v = \sum_{i=1}^{\mu} v_i(t) e_i\]
will be {horizontal if and only if the coefficients $v_i(t)$
solve the system of differential equations
\[ \frac{d}{dt} 
\left( 
\begin{array}{c} 
v_1\\
v_2\\
\ldots\\
v_{\mu} 
\end{array} \right ) =-A(t) \left( \begin{array}{c} v_1\\v_2\\\ldots\\v_{\mu} \end{array} \right ) .\] 

It is of importance to realise that the matrix $A(t)$ will not extend 
holomorphically over $0$. If that were to happen, the coefficients of the 
horizontal sections would be holomorphic at $0$ too, and as a consequence, 
the monodromy would be trivial. So the best one can hope for is a
matrix $A(t)$ with a pole at $0$.\\
We will describe a nice tautological way to do this with a pole of order one.\\

There is no such thing as {\em the} Milnor-fibre, as for each $t \in S^*$
we obtain a different manifold $X_t$. It is somewhat surprising that there 
is an object that might be called the {\em canonical Milnor fibre}: 
one simply takes the pull-back 
$X_{\infty}$ of the fibration $X^* \lra S^*$ over the universal covering 
$\exp:\widetilde{S^*} \lra S^*$, i.e. we form the pull-back diagram
\[
\begin{array}{ccc}
X_{\infty} & \lra &X^*\\
\downarrow&&\downarrow\\
\widetilde{S^*}& \stackrel{\exp}{\lra} &S^*
\end{array},
\]
and we form the groups
\[H_{\Z}:=H^n(X_{\infty},\Z),\;\;H_{\C}:=\C \otimes_{\Z}H_{\Z}=H^n(X_{\infty},\C) .\]
The monodromy will induce an automorphism of $H_{\Z}$ and $H_{\C}$ that we will
denote by $T$.\\
Let us fix an endomorphism $A \in \operatorname{End }(H_{\C})$ with the property that
\[ e^{2\pi i A} T =Id\]
Clearly, if $T$ is a Jordan matrix with eigenvalue $\lambda$, then the choices for 
$A$ are determined by choices of exponents $\alpha$ such that 
\[ e^{-2\pi i\alpha}=\lambda .\]

From the vector space $H_{\C}$ and $A$ we now can construct a specific
extension of the cohomology bundle.
Let $h \mapsto h(t)$ be the isomorphism between $H_{\C}$ and $H^n(X_t,\C)$
induced from the maps
\[ X_t \hookrightarrow X_{\infty}, \;\;\;t \in \widetilde{S^*}\]
Then 
\[ s[h]:=t^A h(t) \]
can be considered as a {\em single valued section} of $\mathcal{H}_{S^*}$:
if we replace $t$ by $e^{2\pi i \varphi} t$ and let $\phi$ increase from $0$ to 
$1$,
we obtain
\[ t^{A}t^{2\pi i A} (T(h))(t)= t^A h(t). \]
Now let 
\[ \mathcal{L} \subset j_*(\mathcal{H}_S) \]
be the $\mathcal{O}_S$-module generated by these $s[h],\;h \in H_{\C}$.
It is called the {\em canonical extension} of $\mathcal{H}_{S^*}$ associated 
to $A$. Note that 
\[ t\frac{d}{dt} t^A =A t^A\]
so that on $\mathcal{L}$ the connection has tautologically a first order pole:
\[ t \partial_t: \mathcal{L} \to \mathcal{L} .\]
Note that the map
\[H_{\C} \lra \mathcal{L};\;\;\; h \to t^A h(t)\]
induces a natural isomorphism
\[ H_{\C} \stackrel{\approx}{\lra} \mathcal{L}/t\mathcal{L}\]
under which the map $A \in \operatorname{End}(H_{\C})$ is identified with $t\partial_t$, the so-called {\em residue} of the connection.

\subsection{The Brieskorn module}

We denote by
\[ \Omega^p:=\Omega^p_{\C^{n+1},0}\]
the $\mathcal{O}$-module of germs at the origin of holomorphic differential 
$p$-forms on $\C^{n+1}$.
The factor space
\[\mathcal{H}^{(0)} :=\Omega^{n+1}/df\wedge d\Omega^{n-1}\]
is called the {\em Brieskorn module}.
Note that what is divided out here is not an $\mathcal{O}$-submodule, so
the Brieskorn module is {\em not} an $\mathcal{O}$-module.
As we will see, this innocent looking algebraic object contains a wealth 
of information. We will see later that $\mathcal{H}^{(0)}$ can be seen as 
stalk at $0$ of a specific extension of the cohomology bundle $\mathcal{H}_{S^*}$.

It is readily verified that the Brieskorn module $\mathcal{H}^{(0)}$ carries {\em two operations}
\[ t : \mathcal{H}^{(0)} \lra \mathcal{H}^{(0)}, [\omega] \mapsto [f\omega]\]
\[ b: \mathcal{H}^{(0)} \lra \mathcal{H}^{(0)}, [\omega=d\eta] \mapsto [df \wedge \eta]\]
which satisfy the relation
\[ tb-bt=b^2 .\]

We now explain how the Brieskorn module $\mathcal{H}^{(0)}$ is naturally associated 
to the period integrals considered above. First note that the period integrals
\[ I(t)=\int_{\delta(t)} Res_{X_t}\left(\frac{\omega}{f-t}\right)\]
only depend on the {\em class} $[\omega]$ of $\omega \in \mathcal{H}^{(0)}$:
for $\omega = df \wedge d\lambda$ we find
\[ I(t)=\int_{\delta(t)} Res_{X_t}\left(\frac{d(f-t)\wedge d\lambda }{f-t}\right)=\int_{\delta(t)} d\lambda=0\]
as the integral of an exact form on a cycle vanishes by Stokes theorem.
Furthermore, 
\[0=\int_{\Delta(t)} \frac{(f-t)\omega}{f-t},\]
so that indeed  
\[\int_{\delta(t)} Res\left(\frac{f\omega}{f-t}\right)=t\int_{\delta(t)}Res\left(\frac{\omega}{f-t}\right) .\]
Of particular importance is the {\em derivative} of a period integral. By
differentiating under the integral sign and using the fact that one can freely
'move' the cycle $\Delta(t)$ without changing the integral, we see that
\[\frac{d}{dt} \int_{\Delta(t)} \frac{\omega}{f-t}=\int_{\Delta(t)} \frac{\omega}{(f-t)^2} .\]
But this is no longer an integral of the type we were considering, as now we 
encounter a pole of order two along $X_t$. However, one has the following key 
fact:\\

{\bf \em Proposition:}
\[\frac{d}{dt} \int_{\Delta(t)} \frac{ b\omega}{f-t}= \int_{\Delta(t)} \frac{\omega}{f-t}\]
{\bf proof:} We write $\omega=d\eta$ and have by definition $df\wedge \eta=b\omega$. Then:

\[\frac{d}{dt} \int_{\Delta(t)} \frac{df \wedge \eta}{f-t}=\int_{\Delta(t)} \frac{df \wedge \eta}{(f-t)^2}=\int_{\Delta(t)} \frac{\omega}{f-t}\]
where we used the easy-to-check {\bf \em pole order reduction formula:}
\[ \frac{df \wedge \eta}{(f-t)^2}=\frac{d\eta}{f-t}-d(\frac{\eta}{f-t}) .\]
The integral over the second term over $\Delta$ clearly gives $0$.\hfill $\Diamond$.\\

This proposition shows what the mysterious operation $b$ on the Brieskorn
module really is: it is to be identified with the {\em inverse} of 
differentiation with respect to $t$. One often writes:
\[ \partial_t^{-1}:=b\]
 
The following result of fundamental importance of was conjectured by {\sc Brieskorn (1970)}:\\
 
{\bf \em Theorem ({\sc Sebastiani (1970)}:)}
{\em If $f$ is an isolated singularity with Milnor-number $\mu$, then
$\mathcal{H}^{(0)}$ is a free $\C\{t\}$-module of rank $\mu$.}\\

We will see later that the freeness is a consequence of the {\em positivity lemma}.\\

As a direct consequence, the factor space
\[\mathcal{H}^{(0)}/t\mathcal{H}^{(0)}=\Omega^{n+1}/(df\wedge d\Omega^{n-1}+f\Omega^{n+1})\]
is a $\C$-vector space of dimension $\mu$. Note that the space 
\[\mathcal{H}^{(0)}/\partial_t^{-1}\mathcal{H}^{(0)}=\Omega^{n+1}/df\wedge \Omega^n=:\Omega_f\]

is also a $\C$-vector space of dimension $\mu$, but there is no
simple canonical isomorphism between the two, unless $f$ is weighted 
homogeneous, in which case we learn from the Euler relation
\[\sum w_i x_i \partial_i f= f\]
that the two subspaces that are divided out are {\em equal}.\\

If we choose a basis
\[\omega_1,\omega_2,\ldots,\omega_{\mu} \]
for $\mathcal{H}^{(0)}$ as $\C\{t\}$-module, one can in principle
write out the action of $b=\partial_t^{-1}$ on this basis and
obtain a $\mu \times \mu$-matrix $B(t)$ of homolomorphic function germs 
in $t$ such that

\[ \partial_t^{-1} \left( 
\begin{array}{c}
\omega_1\\
\omega_2\\
\ldots\\
\omega_{\mu}
\end{array}
\right) = B(t)
\left (
\begin{array}{c}
\omega_1\\
\omega_2\\
\ldots\\
\omega_{\mu}
\end{array} \right)\]

For the corresponding period integrals 
\[I_i :=\int_{\delta(t)} Res(\frac{\omega_i}{f-t})\]
over some vanishing cycle $\delta(t)$ we thus find
\[  
\left( 
\begin{array}{c}
I_1\\
I_2\\
\ldots\\
I_{\mu}
\end{array}
\right) =
\frac{d}{dt}( B(t)
\left (
\begin{array}{c}
I_1\\
I_2\\
\ldots\\
I_{\mu}
\end{array}
\right))
\]

which amounts to a linear system of differential equations for the
period integrals
\[
\frac{d}{dt} 
\left( 
\begin{array}{c}
I_1\\
I_2\\
\ldots\\
I_{\mu}
\end{array}
\right) = A(t)\left( 
\begin{array}{c}
I_1\\
I_2\\
\ldots\\
I_{\mu}
\end{array}
\right),
\]
where
\[A(t)=B(t)^{-1}\left(I-\frac{d}{dt}B(t)\right) .\]
As we know from the {\em positivity lemma} that all period integrals have
moderate growth, it follows from the classical theory of systems of 
linear differential equations that the above linear system is in fact
{\em regular singular}: after an appropriate change of basis we can
achieve that $A(t)$ has at most a  pole of order one:
\[ A(t)=\frac{1}{t}A_{-1} +A_0+A_1t+\ldots\]
 and all solutions admit an expansion of the form
\[ \sum_{\alpha,k} C_{\alpha,k} t^{\alpha}(\log t)^k.\]

{\em Exercises:}\\
\begin{enumerate}
\item Verify that 
\[(t\partial_t-\alpha)^{k} t^{\alpha}(\log t)^k \neq 0,\;\;[(t\partial_t-\alpha)^{k+1} t^{\alpha}(\log t)^k=0,\]
\item Solve the differential equation $t^2\partial_t f =1$ on a punctured disc
$S^*$.\\
\item Give an explicit description of the Brieskorn lattice for $f=y^2+x^3$.
Use a monomial basis $x^ny^mdxdy$ to describe the elements of $\Omega^{n+1}$.
Indicate the relations between these monomials coming from the
equivalences induced by the relations
\[d f \wedge d x^ay^b\]
Find a basis for $\mathcal{H}^{(0)}$. What is the action of $t$, $\partial_t^{-1}$ on the basis elements.   
\item Verify that $tb-bt=b^2$
\item Verify that from $\partial_t t-t \partial_t =1$ one obtains formally
\[t\partial_t^{-1}-\partial_t^{-1}t=\partial_t^{-2} .\]
\item Use the Euler relation 
\[\sum_{i=0}^n w_i x_i\partial_i f =f\]
to show that the two sub-spaces 
\[df\wedge \Omega^{n},\;\;\; \textup{and}\;\;\;df \wedge d\Omega^{n-1}+f\Omega^{n+1}\]
of $\Omega^{n+1}$ are equal if $f$ is weighted homogeneous.

\end{enumerate}

\section{Lecture 4, thursday, august 10, 2017}
In the last lecture we have seen how the Brieskorn module $\mathcal{H}^{(0)}$
leads to an algebraic approach to determine the (complexified) monodromy
transformation. However, the construction was rather ad hoc and today we will
see how the Brieskorn module fits into a standard approach using the relative
deRham complex of a good representative $f:X \lra S$ of our singularity. The
most natural interpretation however is in terms of the {\em Gau\ss -Manin
system}, which corresponds to the direct image in the category of 
$\mathcal{D}$-modules.\\

\subsection{The deRham and relative deRham complexes}
The Poincar\'e lemma states that if  $Z$ is a smooth complex manifold of
dimension $n$, then the constant sheaf $\C_Z$ is resolved by the de 
Rham-complex on $Z$:
\[\C_Z = [\mathcal{O}_Z \stackrel{d}{\lra} \Omega_Z^1 \stackrel{d}{\lra} \ldots\stackrel{d}{\lra} \Omega^{n}_Z] \]
It follows that the cohomology of $Z$ is equal to the {\em hypercohomology} of the
de Rham complex:
\[H^i(Z,\C)=\H^i(\C_Z)=\mathbb{H}^i(\Omega_Z^{\bullet}).\]
In the special case that $Z$ is a Stein space, it follows from Cartan theorem B
that the higher cohomology $H^i(\Omega_Z^j)$ vanishes for $i>0$, so in that case
the cohomology of $Z$ can be described by the complex of global holomorphic differential forms on $Z$:
\[H^i(Z,\C)=H^i(\Gamma(Z,\Omega_Z^{\bullet})) .\]
This applies in particular to the total space $X$ of a good representative and
also to the Milnor fibres $X_t$.\\

If $f: Z \lra T$ is a smooth submersion between complex manifolds of fibre dimensoin $n$, one can use the relative de Rham complex $(\Omega^{\bullet}_{Z/T}, d)$
with terms 
\[\Omega^p_{Z/T} :=\Omega_Z^p/df\wedge \Omega^{p-1}_Z,\]
and differential induced by $d$ to describe the cohomology of the fibres of the map $f$. The relative Poincar\'e lemma states that the relative de Rham complex
resolves the sheaf $f^{-1}\mathcal{O}_T$ of holomorphic functions that are constant along the fibres:
\[ f^{-1}(\mathcal{O}_T)=[\mathcal{O}_Z \lra \Omega^1_{Z/T} \lra \ldots \lra \Omega^n_{Z/T}] .\]
Taking the $i$-th direct image we obtain a description of the
cohomology bundle as the hyperdirect image of relative de Rham:
\[\mathcal{O}_T \otimes R^if_* \C_Z= \mathbb{R}^if_*(\Omega_{Z/T}^{\bullet}) .\] 
In case $Z \lra T$ is a Stein-mapping, the higher direct image sheaves of 
$\Omega_{Z/T}^{p}$ vanish, and we are left with
\[\mathcal{O}_T \otimes R^if_* \C_Z= \mathcal{H}^i(f_*(\Omega_{Z/T}^{\bullet})),\]
i.e. we can describe the cohomology bundle using the cohomology of the global
relative de Rham complex.\\

Let us see what remains from this picture if we apply it to the
map $f:X \lra S$, which is non-submersive at $0$.\\

{\bf \em Theorem ({\sc Brieskorn} (1970)):} For a good representative $f:X \lra S$ of an isolated singularity  we set 
\[\mathcal{H}^p(X/S):=\mathbb{R}^pf_*(\Omega_{X/S}^{\bullet})=\mathcal{H}^p(f_*(\Omega_{X/S}^{\bullet})).\]
Then
\begin{enumerate}
\item $\mathcal{H}^p(X/S)$ is $\mathcal{O}_S$-coherent.\\
\item $\mathcal{H}^p(X/S)|_{S^*}=\mathcal{H}^p_{S^*} =\mathcal{O}_{S^*}\otimes R^p f_*(\C_X)$\\
\item there is an exact sequence of sheaves on $S$:
\[ 0 \lra R^pf_*\C_X \otimes \mathcal{O}_S \lra \mathcal{H}^p(X/S) \lra f_*(H^p(\Omega^{\bullet}_f)) \lra 0\]
where the complex $\Omega_f^{\bullet}$ has terms
\[\Omega_f^p:=\Omega^p/df\wedge\Omega^{p-1}\]
and differential induced by $d$.
\end{enumerate}

The deepest fact here is the appearance of {\em coherence} of the cohomology
in this non-proper situation. It is an application of an important functional 
analytic principle that goes back to {\sc H. Cartan} and {\sc L. Schwarz}: 
on the one hand, shrinking of $X$ does not change the cohomology, but on the other hand such a shrinking is a compact operator, which implies the coherence (Theorem of {\sc Forster-Knorr}).
The second fact is obvious from the previous discussion, and shows that
$\mathcal{H}^p(X/S)$ can be seen as a specific extension of the cohomology 
bundle. The third point exhibits the stalk at $0$ of the sheaf in the middle as the cohomology of a very concrete complex of $\mathcal{O}$-modules.\\  
 
Given this theorem of {\sc Brieskorn}, {\sc Malgrange} has given a purely algebraic
proof that the cohomology $H^p(\Omega^{\bullet}_f)$ is non-zero
only for $p=0$ and $p=n$. This gives an algebraic proof that the
groups $H^p(X_t,\C)$ are only non-vanishing for $p=0$ and $p=n$.\\

The argument runs as follows: take $\omega \in \Omega^p$, representing
a class $[\omega] \in H^p:=H^p(\Omega^{\bullet}_f)$. So we have
\[d \omega=df \wedge \eta\]
for some $\eta \in \Omega^p$. So 
\[ df \wedge d \eta =-d(df\wedge \eta)=-d(d\omega)=0 .\]

For an isolated singularity $f$, the partial derivatives
form a regular sequence and as a consequence, the corresponding
Koszul-complex resolves the Milnor-algebra $\mathcal{O}/J_f$. 
This means that the complex $(\Omega^{\bullet},df \wedge)$ 
is exact in degrees $<n+1$, we may 
conclude from this that $d \eta = df \wedge \lambda$, 
for some $\lambda \in \Omega^{p}$ (This step fails if $p=n$). In this
way we obtain a well defined operation
\[\partial_t :H^p \lra H^p,\;\;\;[\omega]\mapsto [\eta] .\]
If we let 
\[ t: H^p \lra H^p,\;\;\;[\omega] \mapsto [f\omega],\]
we obtain the structure of a module over the ring
\[ \mathcal{D}:=\C\{t\}[\partial_t] . \]
But the operation $\partial_t$ is invertible on $H^p$: given
$[\eta] \in H^p$, the form $df \wedge \eta$ is closed:
\[d(df \wedge \eta)=-df \wedge d\eta=-df \wedge df \wedge...=0 .\] 
If $p \neq 0$ the Poincar\'e-lemma implies that $df \wedge \eta=d\omega$
for some $\omega$. Then $[\omega]$ maps to $[\eta]$ under $\partial_t$.
But this can happen only if $H^p=0$.\\ 

The only interesting group is the  $n$-cohomology $H:=H^n(\Omega_f^{\bullet})$.
It  maps into
\[H':=\Omega^n/df\wedge\Omega^{n-1}+d\Omega^{n-1} , \]
which in turn maps to the Brieskorn module
\[H'' :=\Omega^{n+1}/df\wedge d\Omega^{n+1}=\mathcal{H}^{(0)}\]
via the map
\[ \eta \mapsto df \wedge \eta .\]

\subsection{The Gauss-Manin system}
The modern perspective on these matters is via the theory of 
$\mathcal{D}$-modules, that was developed by {\sc Z. Mebkhout}, 
{\sc M. Kashiwara} and others. We now give a very rough sketch 
of this important theory.\\

On any complex manifold $Z$ there is a
sheaf $\mathcal{D}_Z$ of differential operators on $Z$. This 
is a sheaf of non-commutative rings, locally generated by 
$\mathcal{O}_Z$ and the sheaf of vector fields $\Theta_Z$. 
A $\mathcal{D}_Z$-module is a a sheaf on which $\mathcal{D}_Z$
acts and this notion generalises the notion of vector bundle with 
connection on $Z$. There is a so-called de Rham-functor 
$dR$, which maps a $\mathcal{D}_Z$-module to its {\em de Rham complex}
\[dR(\mathcal{M})=[ \mathcal{M} \stackrel{}{\lra}\mathcal{M}\otimes \Omega_Z^1 \stackrel{}{\lra} \mathcal{M} \otimes \Omega_Z^2 \stackrel{}{\lra}\ldots ] \]
It induces a functor from 
\[D^b_{hol}(Z)  \stackrel{dR}{\lra} D^b_c(Z)\] 
from the derived category of $\mathcal{D}_Z$-modules to the derived category
of constructible sheaves on $Z$.
For example, the $\mathcal{D}_Z$-module $\mathcal{O}_Z$ is mapped to the
ordinary de Rham complex $(\Omega_Z^{\bullet},d)$ of $Z$, which by the Poincar\'e-lemma is quasi-isomorphic to $\C_Z$, the constant sheaf on $Z$.
What is more, there is the formalism of the six operations, which
commute with $dR$. In particular, for a map $f: Z \lra T$ one can form the 
direct image $\int^{\bullet} \mathcal{O}_Z$ which is (a complex of) 
$\mathcal{D}_T$-module(s), which via the de Rham-functor should correspond
to the direct image $R^{\bullet}f_*\C_Z$ of the constant sheaf.\\

Without going into any further detail, the lesson is that there should exist
a natural $\mathcal{D}_S$-module associated to a good representative 
$f:X \lra S$ of an isolated singularity. It is called the {\em Gauss-Manin system} of $f$ and has a very simple concrete description that we will give now.\\

{\bf Definition:} {\em The Gauss-Manin-system of $f$} is defined 
\[\mathcal{G} :=\mathcal{H}^{n+1}(\Omega^{\bullet},\underline{d}) .\]

Let us spell this out in detail. We denote, as before, by 
$\Omega^{\bullet}$ the $\mathcal{O}$-module of germs of holomorphic 
diferential forms on $(\C^{n+1},0)$. We consider an abstract symbol
$D$ and form $\Omega^{\bullet}[D]$, whose elements are polynomials
\[ \sum \omega_k D^k\]
in the symbol $D$, with coefficients from $\Omega^{\bullet}$.
On $\Omega^{\bullet}[D]$ there is defined a differential $\underline{d}$,
which is defined on a basis element $\omega D^k$ as follows: 
\[ \underline{d} (\omega D^k ) = d\omega D^k +df\wedge \omega D^{k+1} .\]
Then one has $\underline{d}\underline{d}=0$, which follows from the
anti-commutation of the two operations
\[d:\Omega^p \lra \Omega^{p+1},\;\;\;df \wedge: \Omega^p \lra \Omega^{p+1} .\]
One can depict $\Omega^{\bullet}[D]$ as a double complex of a particular
shape, that is indicated in the following picture.
\[ 
\begin{array}{cccccccccccccccc}
\ldots&\ldots&\ldots&\ldots&\ldots&\ldots&\ldots&\ldots&\ldots&\ldots&\ldots&\ldots&\ldots&\ldots&\ldots\\
&&\uparrow&&\uparrow&&\uparrow&&&&\uparrow&&\uparrow&&\uparrow&\\
0 &\lra &\mathcal{O} &\stackrel{d}{\lra}&\Omega^1&\lra&\ldots&\lra&\Omega^n&\stackrel{d}{\lra} & \Omega^{n+1}&\lra& 0&\lra&0\\
&&\uparrow&&\uparrow&&\uparrow&&&&\uparrow&&\uparrow&&\uparrow&\\
0&\lra&0& \lra &\mathcal{O} &\stackrel{d}{\lra}&\Omega^1&\lra&\ldots&\lra&\Omega^n&\stackrel{d}{\lra}&\Omega^{n+1}&\lra& 0\\
&&&&\uparrow&&\uparrow&&\uparrow&&\uparrow&&\uparrow&&\uparrow&\\
0&\lra&0&\lra&0&\lra &\mathcal{O} &\stackrel{d}{\lra}&\Omega^1&\lra&\ldots&\lra&\Omega^n&\stackrel{d}{\lra}&\Omega^{n+1}\\
\end{array}
\]
The horizontal arrows are induced by the exterior derivative $d$, the vertical maps are induced by $df \wedge$.
The $k$-th row comes with a power $D^k$, which is suppressed from the diagram.\\

On the complex $(\Omega^{\bullet}[D],\underline{d})$ one can define two
further operations
\[ \partial_t (\omega D^k)=\omega D^{k+1} ,\]
\[ t\cdot (\omega D^k) = f\omega D^k-k\omega D^{k-1} .\]
It is a pleasant exercise to check that these operations commute
with $\underline{d}$ and satisfy
\[ \partial_t t -t \partial_t =1 , \]
so that the cohomology groups of $(\Omega^{\bullet}[D],\underline{d})$ 
are $\mathcal{D}:=\C\{t\}[\partial_t]$-modules.\\

This is the result of taking the direct image in the category of $\mathcal{D}$-modules and it all looks all rather mysterious. But in fact there is a clear
interpretation of the terms in the complex, if we keep in mind the
following correspondence:
\[ \omega D^k \leftrightarrow  Res\left( \frac{k!\omega}{(f-t)^{k+1}}\right)\]
so that we are dealing with an abstract version of general period integrals
\[\frac{1}{2\pi i}\int_{\Delta(t)} \frac{k!\omega }{(f-t)^{k+1}}\]
From this one can reverse engineer the operations $\underline{d}$, $t$ and
$\partial_t$ on $\Omega^{\bullet}[D]$.\\

{\bf \em Proposition:} The cohomology groups $\mathcal{H}^p(\Omega^{\bullet}[D],\underline{d})$
are non-zero only for $p=1$ and $p=n+1$.\\

The interesting cohomology group is the one in degree $n+1$.\\

\subsection{The Brieskorn module and the Hodge-Filtration}
Let us look at the term $\Omega^{n+1}$ sitting at the lower right
corner of the complex $(\Omega^{\bullet}[D]),\underline{d})$ and let us
try to understand the kernel of the map
\[ \Omega^{n+1} \lra \mathcal{G} .\]

{\bf \em Proposition:} The kernel of the above map is exactly
\[ df \wedge d \Omega^{n-1}\]

{\bf proof:} If the element 
\[ \underline{\eta}:=\sum_{k=0}^N \eta_k D^k \in \Omega^{n}[D]\]
maps under to $\omega$ under $\underline{d}$, we get the equations
\[\omega=d\eta_0,\;\;\;d \eta_1+df \wedge \eta_0=0,d \eta_2+df\wedge \eta_1=0,\ldots,  df \wedge \eta_N=0. \]
From the last equation and the exactness of $df \wedge$ we get 
\[ \eta_N =df \wedge \lambda_N\]
When we substitue this in the penultimate equation we get
\[ d (df \wedge \lambda_N)+df \wedge \eta_{N-1}=0,\]
which amounts to 
\[df \wedge (-d\lambda_N+\eta_{N-1})=0\]
or
\[ -d \lambda_N+\eta_{N-1} =df \wedge \lambda_{N-1} .\]
Continuing in this way we see that  there
exist $\lambda_k$ such that
\[ \eta_k=d\lambda_{k+1} +df \wedge \lambda_k .\]
In particular, it follows
\[\omega=d \eta_0 \in df \wedge d\Omega^{n-1} \]

\hfill $\Diamond$.

{\bf \em Corollary:} the image of $\Omega^{n+1}$ in $\mathcal{G}$ is
isomorphic to the Brieskorn module $\mathcal{H}^{(0)}$:
\[ \mathcal{H}^{(0)}  \subset \mathcal{G} \]

If we apply $\partial_t$ to $\mathcal{H}^{(0)}$ we end up 
in a larger space
\[ \mathcal{H}^{(1)}:=\partial_t\mathcal{H}^{(0)} \subset \mathcal{G},\]
which is exactly the image of the sub-complex
\[
\begin{array}{ccc}
 \Omega^{n+1}D & \lra & 0\\
\uparrow df \wedge& &\uparrow\\
\Omega^n &\stackrel{d}{\lra}&\Omega\\
\end{array}
\]
in $\mathcal{G}$.
In general, the image of the sub-complex 
\[ \oplus\oplus_{q=0}^{k-n+1} \Omega^{n+1-k+q}D^q \] 
obtained from vertical truncation of the above double complex
is exactly
\[ \mathcal{H}^{(k)}:=\partial_t^k \mathcal{H}^{(0)} .\]
Recall that there was a natural operation of $\partial_t^{-1}$
on $\mathcal{H}^{(0)}$. We can extend this to $k \in \Z$ and define
\[ \mathcal{H}^{(-k)}:=\partial_t^{(-k)}\mathcal{H}^{(0)}\]
This is the real significance of the Brieskorn-module: it is a
free $\C\{t\}$-submodule of $\mathcal{G}$ that defines a filtration
\[ \ldots \subset \mathcal{H}^{(-2)} \subset \mathcal{H}^{(-1)} \subset \mathcal{H}^{(0)} \subset \mathcal{H}^{(1)} \subset \mathcal{H}^{(2)} \subset \ldots \subset \mathcal{G}\]
on it. It is commonly called, after a shift
\[\mathcal{F}^{n-k} = \mathcal{H}^{(k)} .\]
the {\em Hodge filtration} on $\mathcal{G}$. 
Note that the quotient
\[\mathcal{H}^{(0)}/\mathcal{H}^{(-1)}=\Omega^{n+1}/df \wedge \Omega^n =:\Omega_f \]
is our Milnor-module of dimension $\mu$ and using the isomorphism 
induced by $\partial_t$ one sees that all the succesive quotients
\[\mathcal{H}^{(k)}/\mathcal{H}^{(k-1)}\]
are isomorphic to  $\Omega_f$.

\subsection{The $V_{\bullet}$-filtration}

The $\mathcal{D}$-module $\mathcal{G}$ has three important properties. First, 
it is $\mathcal{D}$-coherent, second it regular singular and third, the
operator
\[ \partial_t :\mathcal{G} \lra \mathcal{G}\]
is an isomorphism. These properties all follow more or less directly from the
results for the Brieskorn-module. The last property is sometimes referred to
as the {\em micro-local} nature of the Gauss-Manin system.
 

The generalised eigenspaces of the operator
\[t\partial_t -\alpha:\mathcal{G} \lra \mathcal{G}\]
are denoted by
\[C_{\alpha}:=\{ m \in \mathcal{G}\;|\;\exists N\;\;\; (t\partial_t-\alpha)^N m=0\} .\]
Furthermore, we denote by
\[ V_{\alpha} \subset \mathcal{G},\]
the $\C\{t\}$-module generated by all $C_{\beta}$ with $\beta \ge \alpha$ and
similarly
\[ V_{>\alpha} \subset V_{\alpha}\]
generated by $C_{\beta}$ with $\beta >\alpha$.
So \[Gr^V_{\alpha}=V_{\alpha}/V_{>\alpha} \approx C_{\alpha} .\]

The positivity lemma of {\sc Malgrange} implies
\[ \mathcal{H}^{(0)} \subset V_{-1}\]
and it is a very non-trivial fact that can be proved using the self-duality
of the Gau{\ss}-Manin system that one also has
\[ V_n \subset \mathcal{H}^{(0)}\]

The spaces $\mathcal{F}^p=\mathcal{H}^{(n-p)}$ induce a filtration on
the spaces $C_{\alpha}$ by putting
\[ F^pC_{\alpha}:=(\mathcal{F}^p \cap V_{\alpha}+V_{>\alpha})/V_{\alpha} .\]

\subsection{The mixed Hodge Structure}
It was shown by Steenbrink that there is a natural mixed
Hodge structure on the cohomology of the canonical Milnor fibre $X_{\infty}$
of an isolated hypersurface singularity. Recall that $X_{\infty}$ is simply 
the pull-back 
$X_{\infty}$ of the fibration $X^* \lra S^*$ over the universal covering 
$\widetilde{S^*} \lra S^*$, i.e. we form the pull-back diagram
\[
\begin{array}{ccc}
X_{\infty} & \lra &X^*\\
\downarrow&&\downarrow\\
\widetilde{S^*}& \lra &S^*
\end{array} .
\]
As before, we put $H_{\Z}:=H^n(X_{\infty},\Z)$ and on it we have the monodromy 
transformation $T:H_{\Z} \lra H_{\Z}$.
We write $T=T_s T_u=T_uT_s$, where $T_s$ is the semi-simple and $T_u$ the 
unipotent part of the monodromy. The monodromy logarithm is
\[ N:=\log T_u: H_{\Q} \lra  H_{\Q},\;\;\;H_{\Q}=\Q \otimes_{\Z} H_{\Z} ,\]
and determines a unique {\em monodromy weight filtration}
\[ M_0 \subset M_1 \subset \ldots M_{2k-1} \subset M_{2k}\]
characterised by the property that $N(M_k) \subset M_{k-2}$ and
\[ N^k: Gr^M_{n+k} \stackrel{\approx}{\lra} Gr^M_{n-k} .\]

Furthermore, if we denote the generalised eigenspace for eigenvalue
$\lambda=\exp(2\pi i \alpha)$, then
\[ H^{n}(X_{\infty},\C)_{\lambda}=C_{\alpha},\]
so we can identify
\[ H^n(X_{\infty},\C)=\bigoplus_{-1 < \alpha \le 0} C_{\alpha},\]
and thus define
\[ F^p H^n(X_{\infty},\C)=\bigoplus_{-1 < \alpha \le 0 }F^p C_{\alpha} .\]

{\bf \em Theorem ({\sc Steenbrink, Scherk}, 1985):} 
With $W_{\bullet}:=M_{\bullet}$, the triple
\[(H_{\Z}, W_{\bullet}, F^{\bullet})\]
is a mixed Hodge structure.\\

This statement means that the filtration $F^{\bullet}$ induces
on all $Gr^W_kH_{\Q}$ a pure Hodge structure of weight $k$; a
$\Q$-Hodge structure of weight $k$ is a $\Q$-vector space $H_{\Q}$
with a filtration $F^{\bullet}$ on $H_{\C}=\C \otimes H_{\Q}$, such that
\[  F^p \oplus \overline{F^{k-p+1}} =H_{\C} .\]  
Equivalently, we have a {\em Hodge decomposition}
\[ H_{\C}=\oplus_{p+q=k} H^{p,q}\;\;\;H^{p,q}:=F^p \cap \overline{F^{q}}=\overline{H^{q,p}},\]
from which one can reconstruct the Hodge filtration as
\[ F^p=\oplus_{p' \ge p} H^{p',k-p'} .\]

The concept of mixed Hodge structure was introduced by {\sc Deligne}, who
showed that each cohomology group $H^k(Z)$ of a possibly singular 
quasi-projective variety $Z$ carries a natural structure of a mixed Hodge
structure. Each algebraic map between quasi-projective varieties induces
a map of mixed Hodge structures. The power of the theory lies in the fact 
that such maps are {\em strictly compatible} with respect to weight and Hodge
filtration: taking $Gr^W_k$ or $Gr_F^p$ are {\em exact functors}.\\

The mixed Hodge structure on $H^n(X_{\infty},\Z)$ defined above is of
quite a different nature as the mixed Hodge structures defined by
{\sc Deligne}; they are called {\em limit or asymptotic} mixed Hodge 
structures.\\

Note that the above definition of the components of the mixed Hodge 
structure attached to a singularity is completely local, but all proofs 
of the Hodge properties use globalisation.\\ 

\subsection{The spectrum}

We now can give a general definition of the spectrum of a singularity:\\

{\bf Definition:}

\[ Sp(f)=s\sum_{-1 < \alpha \le 0}\sum_{p} \dim_{\C} Gr_F^{n-p} C_{\alpha} s^{\alpha+p}\]

So the eigenvectors which belong to the smallest Hodge space $F^n$ are
lifted to exponent in $(0,1]$, those in the next Hodge space $F^{n-1}/F^n$
lift to exponents in $(1,2]$, etc.\\

There is an alternative definition using the $V_{\bullet}$-filtration on the
Milnor module
\[\Omega_f=\Omega^{n+1}/df \wedge\Omega^n=\mathcal{H}^{(0)}/\mathcal{H}^{(-1)}\]

We set
\[ V_{\alpha}\Omega_f:=(V_{\alpha}\mathcal{H}^{(0)}+\mathcal{H}^{(-1)})/\mathcal{H}^{(-1)}\]
It is an exercise to show that there is an isomorphism
\[ Gr_F^{n-p}C_{\alpha} \stackrel{\partial_t^{-p}}{\lra} Gr^V_{\alpha+p} \Omega_f,\]
so that one can alternatively define
\[Sp(f)=s\sum_{\alpha} \dim_{\C} Gr^V_{\alpha}(\Omega_f) s^{\alpha} .\]

The spectrum of a singularity probably hides many deep secrets. A particular nice one is the 
following, that was inspired by an analysis of the genus one partition function
of the Frobenius manifold defined by the singularty $f$:\\

{\bf Variance Conjecture} ({\sc Hertling}):

\[ \frac{1}{\mu}\sum_{i=1}^{\mu}\left(\alpha_i-\frac{n+1}{2}\right)^2 \le \frac{\alpha_{\mu}-\alpha_1}{12}\]

For details we refer to the monograph of {\sc Hertling}.\\

{\em Exercises}:\\
\begin{enumerate}

\item Check the details of the argument of Malgrange for the
vanishing of $H^p(\Omega^{\bullet}_f)$. In particular, show that
a $\mathcal{D}$-module $H$ on which $\partial_t$ is invertible,
and which is finitely generated as on $\C\{t\}$-module has to be
zero.\\

\item Using the correspondence
\[ \omega D^k \leftrightarrow Res\left(\frac{k!\omega}{(f-t)^{k+1}}\right)\]
check that the action of multiplication by $t$ and $\partial_t$
on $\mathcal{G}$ correspond to the multiplication by $t$ and differentiation
with respect to $t$.\\

\item Check that the operations $t, \partial_t$ on the complex $(\Omega^{\bullet}[D],\underline{d})$ are well defined and commute with $\underline{d}$.\\

\item Show that if an element $m \in \mathcal{G}$ is $\C\{t\}$-torsion (i.e.
$t \cdot m=0$, then it it contained in $V_{-1}\mathcal{G}$. Conclude that the Brieskorn lattice is a free $\C\{t\}$-module.\\

\item Analyse the Gauss-Manin system for the $A_1$-singularity. Show that
if the number of variables is even, $\mathcal{G}$ contains $t$-torsion elements.\\

\end{enumerate}

\section{Lecture 5}

We have sketched how from an isolated singularity $f \in \mathcal{O}$ one
can obtain a mixed Hodge structure $(H_{\Z}, W_{\bullet}, F^{\bullet})$. Here
$H_{\Z}=H^n(X_{\infty},\Z)$ is the cohomology of the canonical Milnor fibre,
$W_{\bullet}=M_{\bullet}$ the monodromy weight-filtration and $F^{\bullet}$ the
Hodge-filtration, that can be obtained from the Brieskorn-lattice inside
the Gau{\ss}-Manin system of $f$. We defined the spectrum $sp(f)$, which 
lifts monodromy eigenvalues to exponents whose integer part records the
Hodge space they are contained in. In this way, the Hodge spaces
correspond to integer intervals $(k,k+1]$ in the spectrum.\\

\subsection{Some further results}
For example for a curve singularity the spectrum is contained in $(0,2)$.
One has 
\[ \# (0,1] \cap sp(f) =\delta(f)=\dim F^1,\;\;\# (1,2] \cap sp(f) =\delta(f)-r+1=\dim F^0/F^1\] 
where $r$ denotes the number of local branches of the curve $f=0$.

For a surface singularity $V:=\{f(x,y,z)=0\}$, the spectrum is 
contained in $(0,3)$ and
\[\# (0,1] \cap sp(f)=p_g(V)=\dim F^2\]
Here $p_g(V)$ is the so-called {\em geometrical genus}, which is defined
as the number of independent differential forms of top degree that do not 
extend on a resolution. To be precise, if 
\[ \pi:\widetilde{V} \lra V\]
is a resolution of an $n$-dimensional isolated singularity $V$, and 
\[ E=\pi^{-1}(0)\]
the exceptional divisor, then
\[p_g(V):=\dim H^0(\widetilde{V}\setminus E, \Omega^n_{\widetilde{V}})/H^0(\widetilde{V},\Omega^n_{\widetilde{V}}) .\] 
The singularities with $p_g(V)=0$ are precisely the {\em rational singularities}, so in the surface case these are the  $A-D-E$-singularities. The $A-D-E$-surface singularities can thus be characterised by the property that 
\[sp(f) \subset (1,2).\]

The spectral numbers in the first interval $(0,1]$ are easy to compute
using an {\em embedded resolution.} So if $f:X \lra S$ is a good 
representative of our singularity, we consider $\rho: Y \lra X$, such that
$(f \circ \rho)^{-1}(0)$ is a divisor with normal crossings. Let $E_i$ be
the irreducible components of $\rho^{-1}(0)$. If $\omega \in \Omega^{n+1}$
is a differential form, we put
\[ w(\omega):=\min_{i} mult_{E_i}(\rho^*(\omega)+1)/mult_{E_i}(f\circ \rho) \]
where $mult_{E_i}$ denotes the vanishing order along $E_i$.\\

{\bf \em Theorem ({\sc Varchenko}, 1982):} The spectral numbers in $(0,1] \cap sp(f)$ 
are exactly the numbers $w(\omega)$ which are $\le 1$.\\

These numbers can be defined in much greater generality and also appear 
under the name of {\em jumping coefficients} in the theory of multiplier 
ideals. The smallest spectral number is called classically {\em the complex 
singularity exponent}, but is now often  called the  {\em log canonical 
treshold}.\\

The Newton-diagram of a function $f$ defined a canonical filtration on
$\Omega^{n+1}$ and thus determines an induced filtration on $N_{\bullet}\Omega_f$.
One has the following very useful theorem due to {\sc M. Saito}:\\

{\bf \em Theorem ({\sc M. Saito}, 1988):} If $f$ is {\em non-degenerate with respect to its 
Newton-diagram}, the Newton filtration $N_{\bullet}$ on  $\Omega_f$ coincides 
with the $V_{\bullet}$-filtration, shifted by one:
\[ V_{\alpha} \Omega_f= N_{\alpha+1}\Omega_f .\]

\subsection{Adjacencies}
It is of great beauty and interest to study how complicated critical
points can decompose into simpler ones. One considers holomorphic
function germs $f_{\lambda}(x)=F(x,\lambda) \in \C\{x_0,x_1,\ldots,x_n,\lambda\}$
 that depend on an additional deformation parameter $\lambda$. If for
$\lambda \neq 0$, small enough, the function $f_{\lambda}$ has critical
points $p_1,p_2,\ldots, p_n$ with critical value $0$ and $(f_{\lambda},p_i) \sim (g_i,0)$, we say that $F$ defines an {\em adjacency} between $f$ and $g_1,g_2,\ldots, g_N$ and write $f \rightsquigarrow g_1,g_2,\ldots,g_N$.\\

{\em Examples:}\\
\begin{center}
\includegraphics[height=4cm]{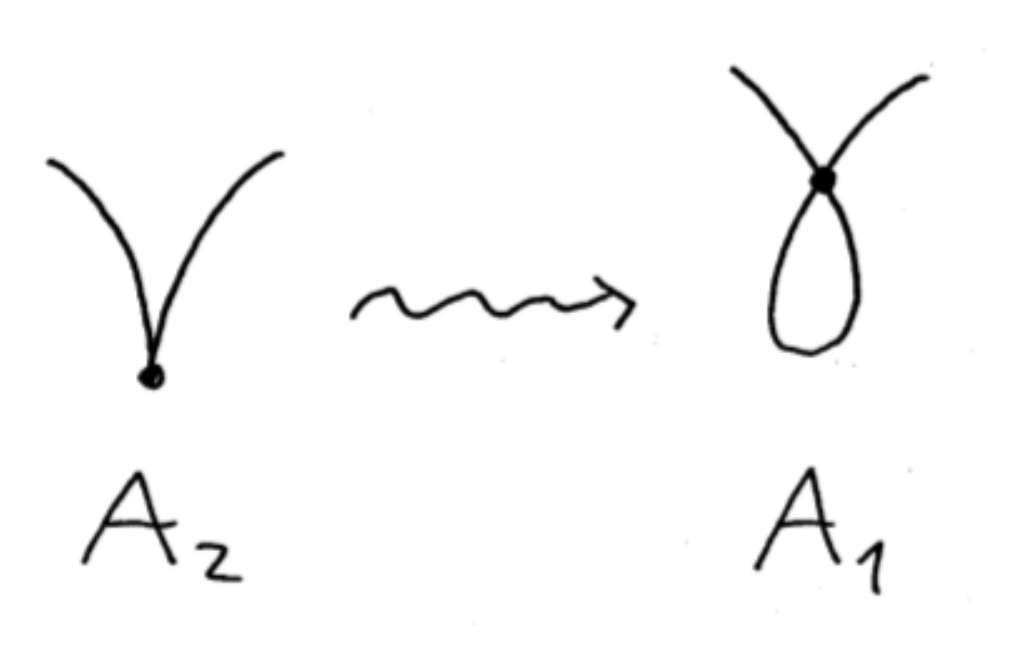}
\includegraphics[height=4cm]{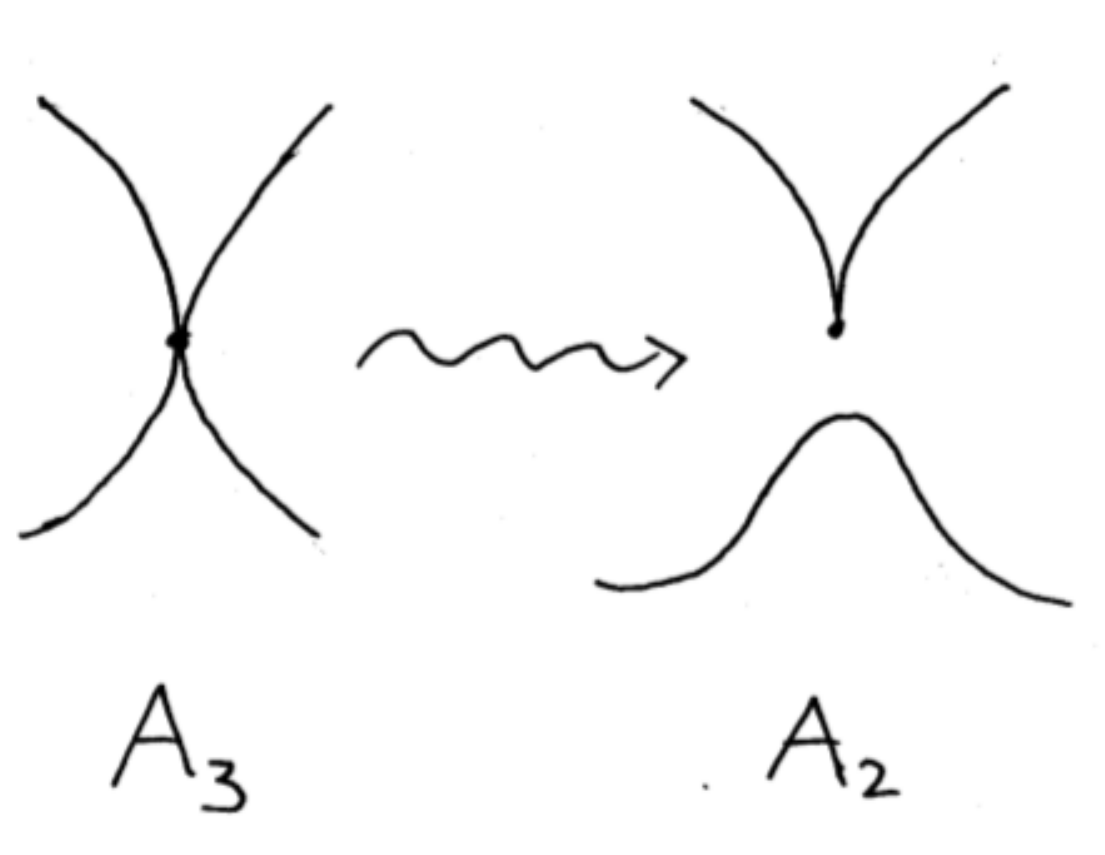}
\includegraphics[height=4cm]{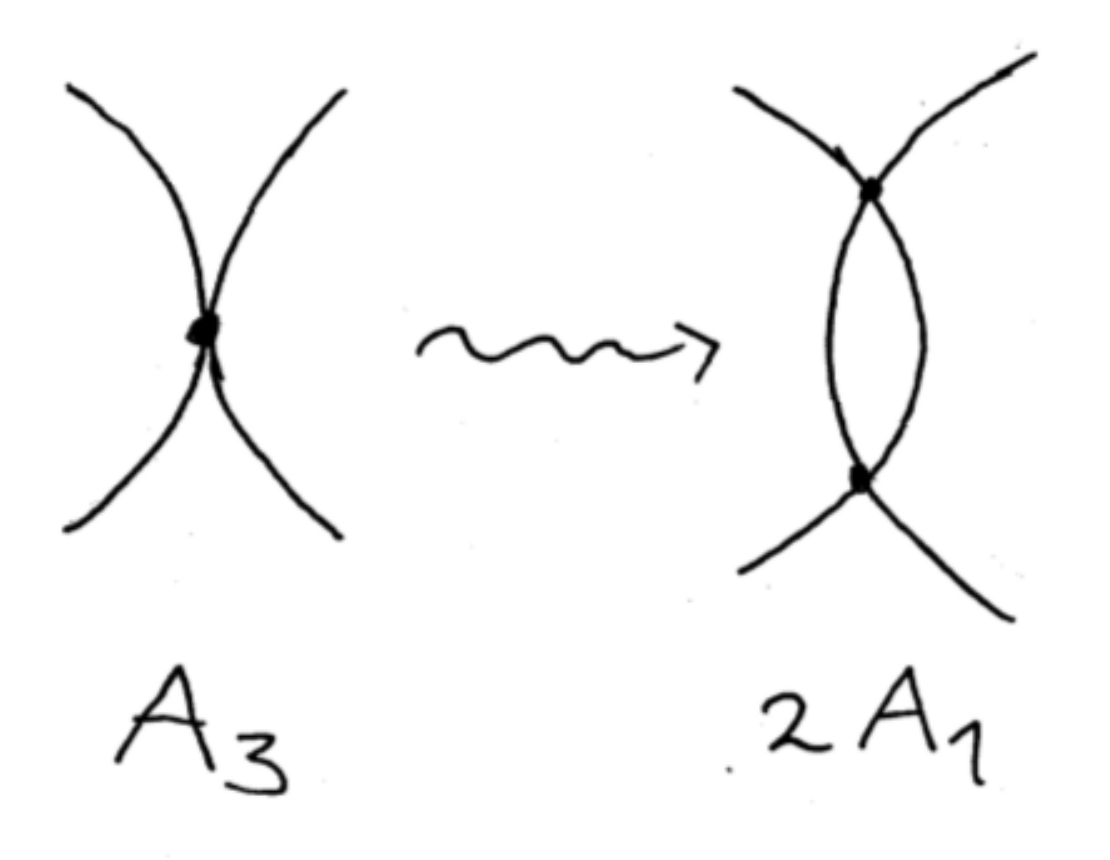}
\includegraphics[height=4cm]{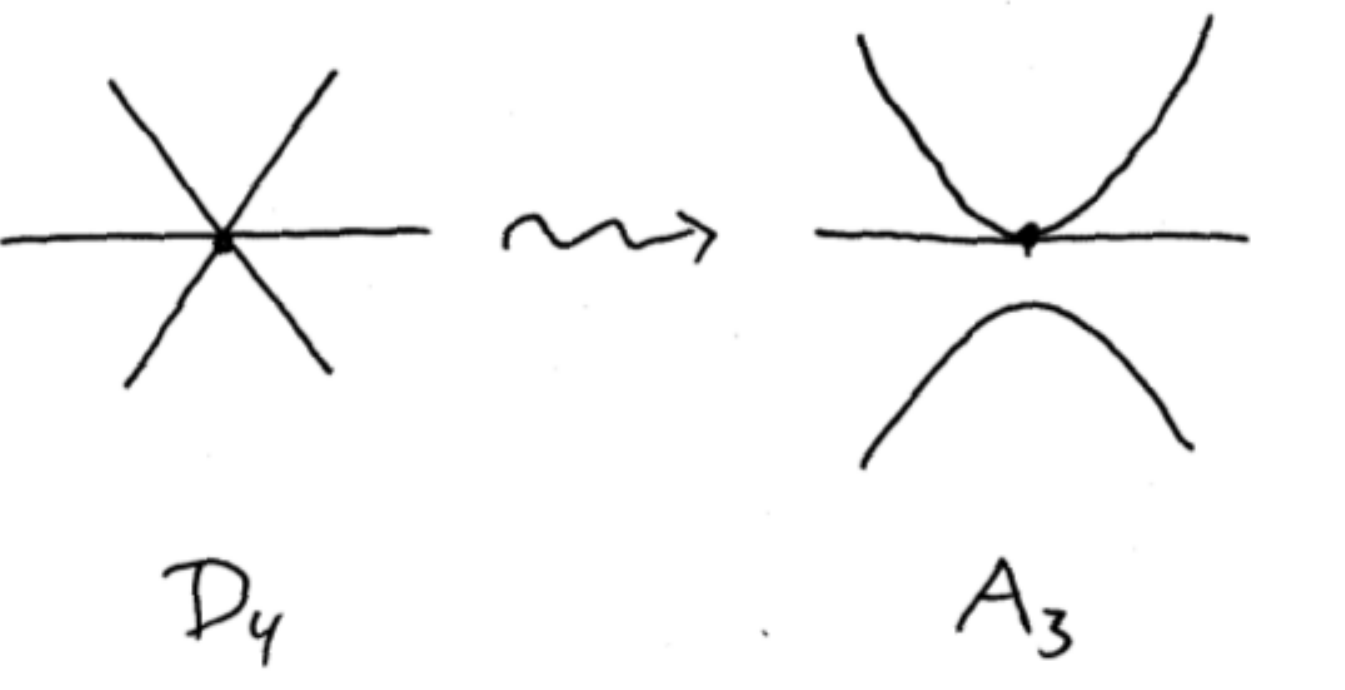}
\includegraphics[height=4cm]{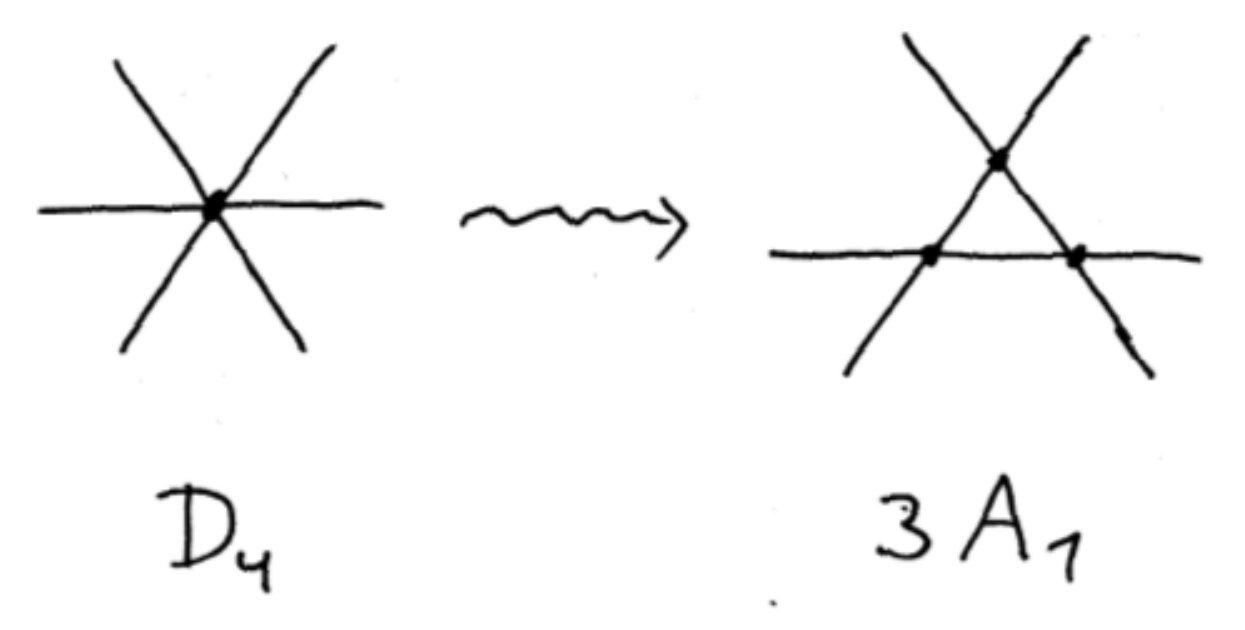}
\includegraphics[height=4cm]{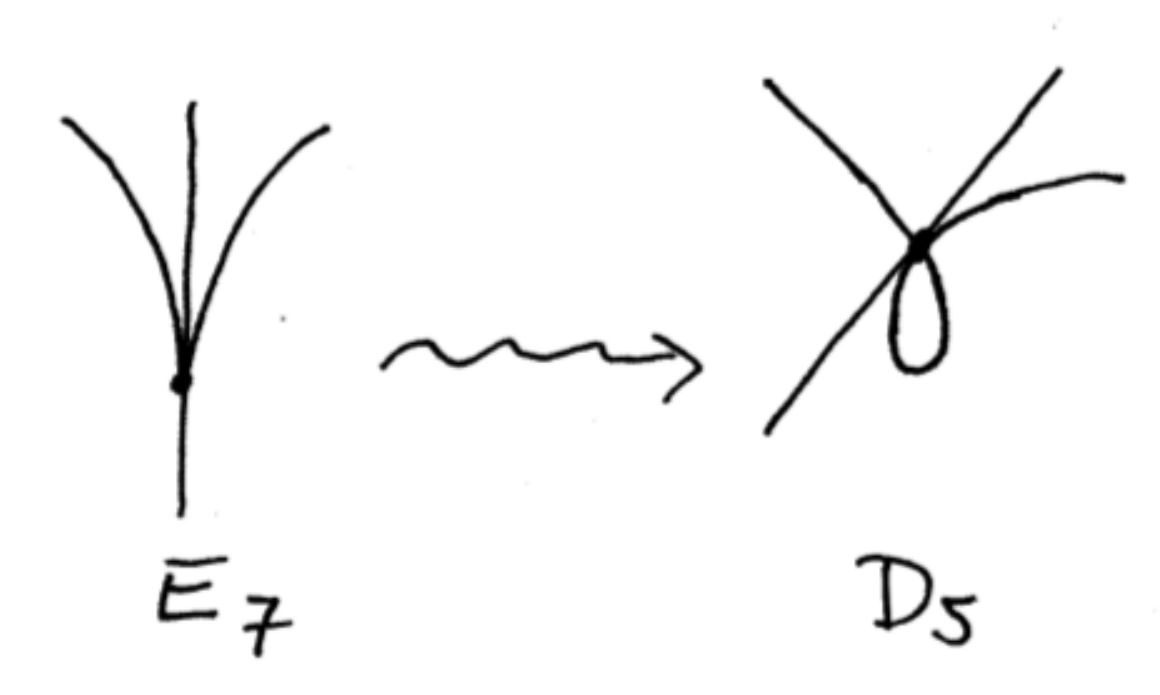}
\includegraphics[height=4cm]{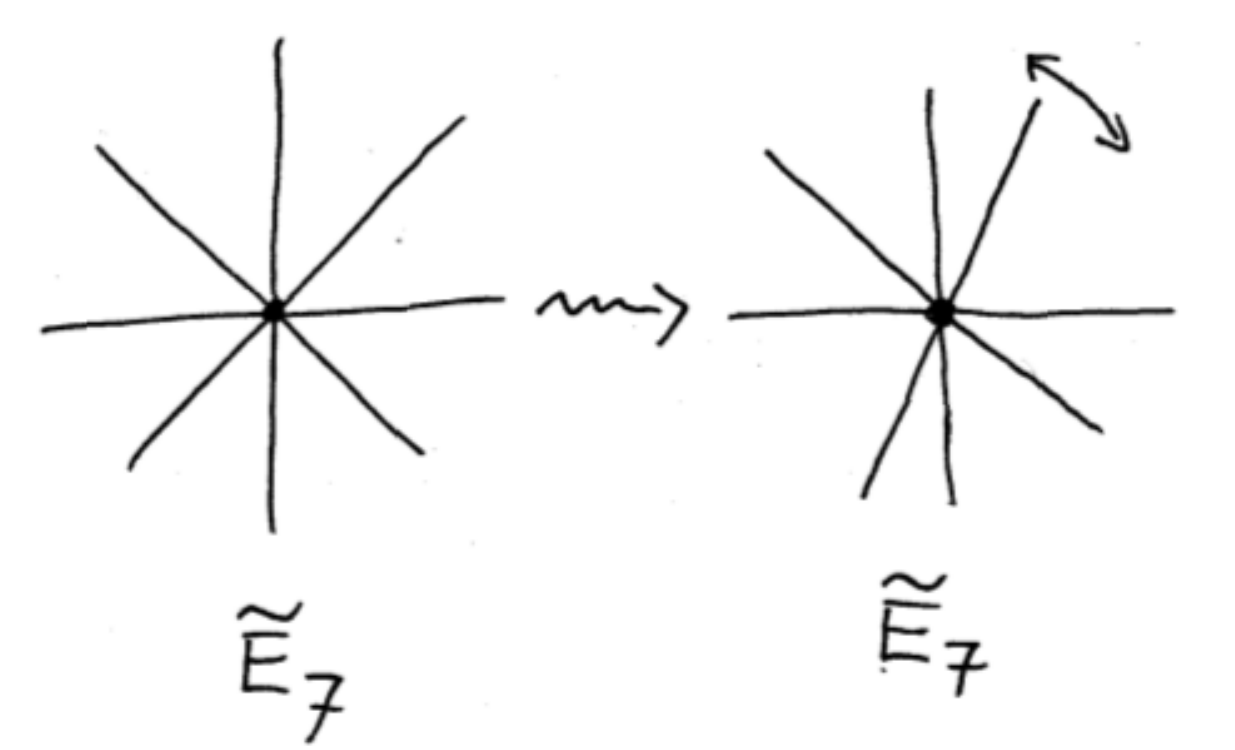}
\end{center}

It is an instructive exercise to write down the formulas that
that describe the above adjacencies. 

\subsection{Semi-continuous invariants}
An invariant $I$ of a singularity is {\em semi-continuous} is
for each adjacency $f \rightsquigarrow g_1+g_2+\ldots+g_N$ one
has
\[ I(f) \ge \sum_{i=1}^N I(g_i)\]

{\em Examples}:
\begin{enumerate}
\item $\mu(f) \ge \sum_{i=1}^{N} \mu(g_i)$\\
\item  $\delta(f) \ge \sum_{i=1}^{N} \delta(g_i)$\\
\item  $p_g(f) \ge \sum_{i=1}^{N} p_g(g_i)$\\
\end{enumerate}
{\sc Arnol'd} conjectured that the quantity
\[ \# (-\infty,a] \cap sp(f)\]
is semi-continuous.\\

{\bf \em Definition:} A subset $S \subset \R$ is called a semi-continuity set if
\[ \# S \cap sp(f)\]
is semi-continuous.\\

{\bf \em Theorem:}\\

{\em (A. Varchenko, 1983):if $f$ is quasi-homogeneous then each open interval $(a,a+1)$
is a semi-continuity set.\\

(J. Steenbrink, 1985): for general $f$ each interval $(a,a+1]$ is semi-continuity set.}\\

{\bf \em Corollary:} In a $\mu$-constant family the spectrum is constant.\\

{\bf proof:} By semi-continuity, the number of exponents( =spectral numbers) in each interval  of length $1$ can only increase under specialisation. As the total number  of exponents is equal to $\mu$ and thus stays constant by assumption, in each interval of lenght $1$ the number of exponents is constant, which can only happen if all exponents stay constant!\\ 

\subsection{The Bruce deformation}
Consider a projective hypersurface
$$Z=\{F(x_0,x_1,\ldots,x_n)=0\}\subset \P^n$$ 
where $F$ is a homogeneous polynomial in $x_0,x_1,\ldots,x_n$. We
can write
\[F=f_d+x_0 f_{d-1}+x_0^2 f_{d-2}+\ldots+x_0^d\]
where $f_k$ is homogeneous of degree $k$ in $x_1,x_2,\ldots,x_n$.
By putting $x_0=\lambda \in \C$, we obtain the equation for an affine 
hypersurface 
\[ X_{\lambda} :=\{ F_{\lambda}=0 \} \subset \C^n, \] 
where
\[ F_{\lambda}:=f_d + \lambda f_{d-1}+\lambda^2 f_{d-2}+\ldots+\lambda^d \in \C[x_1,x_2,\ldots,x_n]\]
Note however that 
\[ F_{\lambda}(\lambda x_1,\ldots,\lambda x_n)=\lambda^d F_1(x_1,\ldots,x_n),\]
which implies that for all $\lambda \neq 0$ these affine hypersurfaces 
$X_{\lambda}$ are isomorphic to each other. As 
\[ X_0 =\{ f_d=0\} \subset \C^n\]
we see that $X_0$ is the affine cone on the projective hypersurface
$Z \cap H$, where $H$ is the hyperplane defined by $x_0=0$.\\ 

If $Z$ has isolated singularities with singular set 
$\Sigma:=\{p_1,p_2,\ldots,p_N\}$ and coordinates are chosen such that
$H \cap \Sigma =\emptyset$, then the homogeneous affine hypersurface
$X_0$ defined by $f_d=0$ defines an isolated singularity, namely the
affine cone over the smooth projective variety $Z \cap H$, whereas
for $\lambda \neq 0$ the affine hypersurface $X_{\lambda}$ has isolated
singularities, defined by $g_1,g_2,\ldots,g_N$
Hence, we obtain an special adjacency
\[ f_d \rightsquigarrow g_1+g_2+\ldots+g_N,\]
that was already used by {\sc J. W.  Bruce}.\\

{\bf \em Corollary:} The spectral bound for projective hypersurfaces holds.\\

The Varchenko semi-continuity theorem implies
\[\# (a,a+1) \cap sp(f_d) \ge \sum \# (a,a+1) \cap sp(g_i)\] 
But $f_d$ and $x_1^d+x_2^d+\ldots+x_n^d$ are connected via a $\mu$-constant
deformation and thus have the same spectrum, so that we have 
\[\# (a,a+1) \cap sp(x_1^d+x_2^d+\ldots+x_n^d) \ge \sum \# (a,a+1) \cap sp(g_i)\] 
As the spectrum at the left hand side is determined by the lattice point in
a cube, we obtain the spectral bound as formulated in lecture 1.

\subsection{Globalisation and variations of (mixed) Hodge structures}

It follows from the fact that any isolated singularity is right equivalent to a polynomial that the Milnor fibration $f:X \lra S$ can be globalised to 
a flat, projective family 
\[ F: Y \lra S\]
and we may assume that the fibre $Y_t=f^{-1}(t)$ over $t$ is smooth for $t \neq 0$ and $Y_0$ has a single isolated singular point, and after restriction we obtain our good
representative $f:X \lra S$.\\
 
We would like to express the idea that the difference between $Y_0$ and $Y_t$ is concentrated at the singularity. This can
be done in great generality with the diagram that can be found in SGA 7.
\[
\begin{array}{ccccccc}
Y_0&\stackrel{i}{\lra}&Y&\stackrel{j}{\longleftarrow} &Y^*&\stackrel{\rho}{\longleftarrow}&Y_{\infty}\\
\downarrow&&\downarrow&&\downarrow&&\downarrow\\
\{0\}&\lra&S&\longleftarrow&S^*&\longleftarrow&\widetilde{S^*}\\ 
\end{array}
\]

The cohomology group $H^n(Y_t)$ of the smooth projective varieties $Y_t$, $t \in S^*$ each carry a pure Hodge structure and form what is called a {\em variation of Hodge structures}. This means that on the cohomology bundle
$\mathcal{H}$ on $S^*$ we are given a filtration by vector bundles on $S^*$
\[ \mathcal{F}^n \subset \mathcal{F}^{n-1} \subset \ldots \subset \mathcal{F}^0=\mathcal{H}\]
satisfying {\em Griffiths transversality}:
\[ \partial_t (\mathcal{F}^k) \subset \mathcal{F}^{k-1}\]
The construction of {\sc W. Schmid} extends these vector bundles to 
$\mathcal{F}^p_e$ on all of $S$. Together with the monodromy filtration $M_{\bullet}$ coming from $N=\log T_u$, one obtains a {\em limiting mixed Hodge structure} on the cohomology groups $H^k(Y_{\infty})$, which were described geometrically using an semi-stable model for the degeneration by J. Steenbrink.
One obtains an exact sequence of mixed Hodge structures
\[ \ldots H^n(Y_0) \lra H^n_{\lim}(Y_{\infty}) \lra H^n(X_{\infty}) \lra H^{n+1}(Y_0) \ldots \]
Here the MHS at the left group is {\sc Deligne}'s MHS on the singular space $Y_0$; in
the middle we have the MHS of {\sc Schmid} and {\sc Steenbrink}. The group at the right is attached
to the singularity, and the MHS of {\sc Scherk} and {\sc Steenbrink} makes it into an 
exact sequence of MHS.\\

For our application we have to dig still deeper. If the general fibre
$Y_t$ is singular itself, we end up with a {\em variation of MHS} with fibres
$H^n(Y_t)$ and one has to consider the problem of how to extend this 
structure over the origin. It turns out there is no direct analogue of the
theorem of {\sc W. Schmid}, rather there is an additional condition leading to the notion of an {\em admissible variation of mixed Hodge structures} introduced
by {\sc J. Steenbrink} and {\sc S. Zucker}, which is satisfied in the geometrical case. The admissibility amounts to existence of a {\em relative weight filtration} $M_{\bullet}$ attached to the monodromy logarithm $N$ and weight filtration $W_{\bullet}$ on $H^n(Y_{\infty})$, characterised by 
\begin{itemize}
\item $N M_k \subset M_{k-2}$
\item $M_{\bullet}$ induced on $Gr^W_k$ the monodromy weight filtration of the
nilpotent transformation $Gr_k N: Gr^W_k \lra Gr^W_k$. 
\end{itemize}
So for each $k$ one obtaines well defined limit mixed Hodge structures
$Gr^W_kH^n(Y_{\infty})$, whose weight filtration comes from the induced monodromy
logarithm.\\

\subsection{The argument}
We consider a flat projective family of varieties with $n$-dimensional
fibres over a base $T$. After {\sc Deligne}, the group $H^i(Y_t)$  
carries for each $t \in T$ a natural mixed Hodge structure.
For each integer $p$ we consder the Euler-characteristic
\[\chi_p(t):=(-1)^n\sum_{i} (-1)^i Gr_F^p H^i(Y_t)\]
It follows from general arguments that this function is {\em construcible}, 
i.e. there exists a stratification for which $\chi_p$ is constant on the 
strata. 
If $h:S \lra T$ maps all points of the punctured disc $S^*$ into one
stratum, and $0$ into an adjacent stratum, then we put
\[\chi_p(h):=(-1)^n \sum_i \dim Gr_F^p \R^i\Phi_h .\]

{\bf \em Proposition:} One has the following {\em jump formula:}
\[\chi_p(h) =\chi_p(t)-\chi_p(s)\]
where $t=h(0)$, $s=h(t')$, $t' \neq 0$.

{\bf proof:} From the triangle of nearby and vanishing cycles
\[\ldots \lra \C_{Y_t} \lra \R\psi_h \lra  \R\phi_h \stackrel{+1}{\lra} \ldots\]
one obtains a long exact sequence
\[\ldots \lra H^k(Y_t) \lra H^k_{lim}(Y_{\infty}) \lra R^k\phi_h \lra H^{k+1}(Y_t) \ldots\]
Here the middle terms is the limit of the variation of mixed Hodge structures
of the family over $h$. Now we use that 
\[\dim Gr_F^p H_{lim}(Y_{\infty}) =\dim Gr_F^P H^p(Y_s)\]
Using this, and taking $\chi_p$ of the above exact sequence then give the
jump formula \hfill $\Diamond$\\

{\bf Definition:} The map $f:Y \lra T$ is called {\em spherical} if
for all maps $h:S \lra T$ as above one has:
\[ Gr_F^p \R^i\phi_h =0\;\;\textup{for}\;\;i \neq n .\]

Examples are families with at most isolated hypersurface singularities,
or more generally, isolated complete intersection singularities. From
the jump formula one gets immediately:\\
 
{\bf Corollary:} If $f:Y \lra T$ is spherical, then the function
\[  t \mapsto \chi_p(t)\]
is upper semi-continuous.\\

Now let $s$ be a generic point of $T$; the difference  between 
$\chi_p(t)-\chi_p(s)$ is also upper semi-continuous and furthermore
can be identified with local contributions coming from the singularities.
 
If we consider an adjacency $f  \rightsquigarrow g_1 +g_2+\ldots g_N$, 
we obtain
\[Gr_F^p H^n(X_{\infty}(f)) \ge \sum_{i=1}^N Gr_F^P H^n(X_{\infty}(g_i))\]
which means
\[ \# (n-p,n-p+1] \cap sp(f) \ge \sum_{i=1}^N \#(n-p,n-p+1] \cap sp(g_i)\]

So all half open intervalls $(k, k+1]$ are semi-continuity sets.

\subsection{Varchenko's trick}
To extend the above argument to arbitrary intervals $(\alpha,\alpha+1]$
one has somehow to shift the intervals. This can be achieved with the
following {\em trick}: The spectrum of
\[ F(x_0,x_1,\ldots,x_n,z)=f(x_0,x_1,\ldots,x_n)+z^m\]
consists of $m$-shifted copies of the spectrum of $f$:
\[Sp(F) =sp(f)+\frac{1}{m} \cup sp(f)+\frac{2}{m} \cup \ldots \cup sp(f)+\frac{m-1}{m} \]
Furthermore, there is an action of $\mu_m$ on the vanishing cohomology of $F$,
induced by 
\[  z \mapsto \zeta z\]
where $\zeta \in \mu_m$ is an $m$-th root of unity. Using now the equivariant
version of the previous set-up gives that the intervals 
$(k+\frac{i}{m},k+1+\frac{i}{m}]$
are also semi-continuity sets. As the spectral numbers are rational, this
clearly implies that all intervals $(\alpha,\alpha+1]$ are semi-continuity
sets.\\

{\em Exercises:}\\

1)  If $f \in \C\{x_0,x_1,\ldots,x_n\}$ is an arbitrary (isolated) singularity, then $F:=f+u^2+v^2$ defines a rational singularity. So $3$-dimensional
 rational singularties are at least as complex as arbitrary curve 
singularities.\\

2) Compute the spectrum of A'Campo's singularity $x^2y^2+x^5+y^5$.\\ 

(i) using the theorem M. Saito's theorem.\\

(ii) using the embedded resolution and Varchenko's theorem.\\

3) Work out the details of {\sc Varchenko}'s trick.\\

We have come to the end of these lectures, which really only pretend to
offer a rough sketch of the ideas involved. I hope the reader will be sufficiently
intrigued to make a more serious study of the original literature. 
I consider myself lucky for having been a witness of a part of the above 
described developments that define a unique {\em golden era of singularity theory}.\\

{\bf Acknowledgement:} First I would like to thank Xavier Gomez-Mont and the
other organisers of the ELGA 3 conference at CIMAT for inviting me to 
give this series of lectures and the audience for asking many questions.
Further thanks to David Mond, Tom Sutherland and Matthias Zach for reading these written notes,
spotting many mistakes and suggesting certain improvements. Of course, I alone
remain responsible for the remaining errors.\\

\section{Literature}
The best way to learn a subject is to study the original papers. 
Only from reading these one can understand the motivations and 
see how the main ideas developed. It is a gratifying aspect of 
the internet era  that it has become rather simple to get acces 
to most of these original papers.\\

{\sc Lecture 1}

\begin{itemize}

\item L. Schl\"afli, {\em On the Distribution of Surfaces of the Third Order into Species, in Reference to the Absence or Presence of Singular Points, and the Reality of Their Lines}, Philosophical Transactions of the Royal Society of London, The Royal Society, {\bf 153}, 193 - 241, (1863).\\
 
\item A. W. Basset, {\em The maximum number of double points on a surface}, Nature, p. 246, (1906).\\

\item J.W. Bruce, C.T.C. Wall, {\em On the classification of cubic surfaces}, 
Journal of the London Mathematical Society, {\bf 19} (2): 245 - 256, (1979).\\

\item J. W. Bruce, {\em  An Upper Bound for the Number of Singularities on a Projective Hypersurface}, Bull. London Math, Soc. {\bf 13}, 47-50, (1981).\\

\item V. I. Arnold, {\em Some problems in singularity Theory},
Proc. Indian Acad. of Schience, Vol {\bf 90}, Nr.1, 1-9, (1981).

\end{itemize}

{\em Lecture 2:}

\begin{itemize}

\item F. Pham, {\em Formules de Picard-Lefschetz g\'en\'eralis\'ees et ramification des integrales}, Bull. Soc. Math. de France {\bf 93}, 333-367, (1965).\\

\item C. H. Clemens, {\em Picard-Lefschetz theorem for families of nonsingular algebraic varieties acquiring ordinary singularities}, Transactions of the American Math. Soc., {\bf 136}, 93-108, (1969).\\

\item P. A. Griffith, {\em Periods of integrals on algebraic manifolds: Summary of main results and discussion of open problems}, Bull. Amer. Math. Soc. {\bf 76},  228-296, (1970).\\

\item N. A'Campo, {\em Sur la monodromie des singularit\'es, isol\'ees d'hypersurfaces complexes}, Inventiones Math. {\bf 20}, 147-169, (1973).\\

\item B. Malgrange,  {\em Letter to the editors}, Inventiones Math. {\bf 20}, 171-172, (1973).\\

\item L\^{e} Dung Trang, {\em Finitude de la monodromie locale des courbes planes}, In: Fonctions de Plusieurs Variables Complexes, S\'eminair F. Norguet, Okt. 1970-D\'ec. 1973, Lectuere Notes in Math. {\bf 409}, Berlin, New York, (1974).\\

\end{itemize}

{\em Lecture 3:}

\begin{itemize}

\item N. Katz, {\em The regularity theorem in algebraic geometry}, Actes, Congr\'es intern. math., Tome 1, 437-443, (1970).\\

\item E. Brieskorn, {\em Die Monodromie der isolierten Singularit\"aten von Hyperfl\"achen}, Manuscripta Math. {\bf 29}, 103-162, (1970).\\

\item  M. Sebastiani, {\em Preuve d'une conjecture de Brieskorn}, Manuscripta Math. {\bf 2}, 301-308, (1970).\\

\item A. Landman, {\em On the Picard-Lefschetz transformation for
algebraic manifolds acquiring general singularities}, Transactions of the American Mathematical Society, Volume {\bf 181}, 89 - 126, (1973).\\

\item B. Malgrange, {\em Int\'egrale asymptotique et monodromie}, Annales scientific de l'ENS, 4e s\'erie, tome 7, no. 3,  405 - 420, (1974).\\

\item M. Saito, {\em On the structure of Brieskorn lattice.} 
Ann. Inst. Fourier (Grenoble) {\bf 39}, no. 1, 27 - 72, (1989).\\ 

\item J. Scherk, {\em On the Monodromy Theorem for Isolated Hypersurfae Singularities}, Inventiones Math. {\bf 58}, 289 - 301, (1980).\\
 
\item M. Schultze, {\em Monodromy of Hypersurface Singularities},
Acta Appl. Math. {\bf 75}, 3 - 13, (2003).\\

\end{itemize}

{\em Lecture 4:}

\begin{itemize}

\item A. Varchenko, {\em The asymptotics of holomorphic forms determine a mixed Hodge structure}, Sov. Math. Dokl. {\bf 22}, 772 - 775, (1980).\\

\item A. Varchenko, {\em On the monodromy operator in vanishing cohomology and the operator of multiplication by $f$ in the local ring}, Sov. Math. Dokl. {\bf 24} 248 - 252, (1981).\\

\item A. Varchenko, {\em Asymptotic Hodge structure in the vanishing cohomology}, Math. USSR. Izv. {\bf 18}, 469 - 512, (1982).\\

\item  J. Scherk, J. Steenbrink, {\em On the Mixed Hodge Structure on the Cohomology of the Milnor Fibre}, Math. Ann. {\bf 271},  641 - 665, (1985).\\

\end{itemize}

{\em Lecture 5:}
\begin{itemize}

\item J. H. M. Steenbrink, {\em Limits of Hodge structures}, Inventiones Math. {\bf 31},  229 - 257, (1976).\\

\item J. H. M. Steenbrink, {\em Mixed Hodge structure on the vanishing
cohomology}, In: Real and Complex Singularities, Oslo 1976, Ed. P. Holme,
Sijthoff-Noordhoff (1977).\\

\item A. Varchenko, {\em On semi-continuity of the spectrum and an upper bound for the number of singular points in projective hypersurfces}, Doklady Ak. Nauk, {\bf 270} (6), 1294 - 1297, (1983).\\

\item J. Steenbrink, S. Zucker, {\em Variations of mixed Hodge structures I},
Inv. Math. {\bf 80}, 489 - 542, (1985).\\

\item J. Steenbrink, {\em Semicontinuity of the singularity spectrum}, Inv. math, {\bf 79}, 557 - 565, (1985).\\

\item M. Saito, {\em Exponents and Newton polyhedra of isolated hypersurface singularities}, Math. Ann. 281, no. 3, 411 - 417, (1988). \\ 

\end{itemize}

Textbooks:\\

\begin{itemize}

\item E. Picard, G. Simart, {\em Th\'eorie des fonctions alg\'ebrique de deux variables ind\'ependantes}, Tome I,  Gauthier-Villars, Paris, (1897).\\

\item S. Lefschetz, {\em L'Analysis situs et la g\'eom\'etrie alg\'ebrique}, Gauthier-Villars, Paris, (1924).\\

\item J. Milnor, {\em Singular points of complex hypersurfaces}, 
Ann. Math. Studies {\bf 61}, Princeton University Press, (1968).\\

\item P. Deligne, N. Katz, {\em SGA7.II: Groupes de Monodomrie en g\'eom\'etrie alg\'ebrique}, Springer Lecture Notes in Mathematics {\bf 340}, (1973).\\

\item F. Pham, {\em Singularit\'es Des Syst\`emes Diff\'erentiels
De Gauss-Manin}, Progress in Mathematics, Vol. {\bf 2}, Birkh\"auser, (1979).\\

\item V.I. Arnold, S.M. Gusein.Zade, A.N. Varchenko, {\em Singularities of differentiable maps, Vol I and II}, Birkh\"auser, (1985).\\

\item E. Brieskorn und H. Kn\"orrer, {\em Ebene Algebraische Kurven}, Birkh\"auser, (1981).\\ 

\item V. Kulikov, {\em Mixed Hodge Structures and Singularities}, Cabridge Tracts in Mathematics {\bf 132}, (1998).\\

\item C. Hertling, {\em Frobenius manifolds and moduli spaces for singularities}, Cambridge Tracts in Mathematics {\bf 151}, (2002).\\

\item C. Peters, J. Steenbrink, {\em Mixed Hodge Structures}, Ergebnisse der Mathematik und ihre Grenzgebiete, Dritte Folge, Vol. {\bf 52}, Springer Verlag, (2008).\\

\end{itemize}

\end{document}